\documentclass [11pt]{article} 
\usepackage[english]{babel}

\usepackage{appendix}
\usepackage{xcolor}
\usepackage{geometry}
\usepackage{amsthm,amsmath,amssymb}
\usepackage{algorithm}
\usepackage{algorithmic}
\usepackage{hyperref}
\usepackage{bbm}
\usepackage{mathrsfs}
\usepackage{graphicx}  
\usepackage{epstopdf}
\usepackage{float}
\usepackage{comment}
\usepackage{url}
\usepackage{tikz}
\usepackage{caption}
\usepackage{graphicx} 
\usepackage{float} 
\usepackage{subfigure}
\usepackage{booktabs}
\usepackage{indentfirst}
\usepackage{authblk}

\newtheorem{theorem}{Theorem}[section]
\newtheorem{lemma}{Lemma}[section]
\newtheorem{Alemma}{Lemma}[section]

\newcommand{\BREAK}{\STATE \textbf{break}}
\newtheorem{corollary}{Corollary}[section]

\newtheorem{Proposition}{Proposition}[section]
\newtheorem{remark}{Remark}[section]

\numberwithin{figure}{section}
\numberwithin{algorithm}{section}

\usepackage{xpatch}
\makeatletter
\AtBeginDocument{\xpatchcmd{\@thm}{\thm@headpunct{.}}{\thm@headpunct{}}{}{}}
\makeatother

\title{\vspace{-0.5in} \bf Service Rate Control in Queues with Abandonments  \vspace{10pt}}
\author{{\normalsize Runhua Wu and  Hayriye Ayhan} \\ \small H. Milton Stewart School of Industrial and Systems Engineering\\
\small Georgia Institute of Technology\\
\small Atlanta, GA 30332-0205, U.S.A.}

\begin{document}
\maketitle

\begin{abstract}
We consider a Markovian single server queue with impatient customers. There is a customer abandonment cost and a holding cost for customers in the system. We consider two versions of the problem. In the first version, customers pay a reward at the time of arrival whereas in the second version, reward is received at the time of service completion. Service rate attains values in a compact set and there is a  cost associated with each service rate.
Under these assumptions, our objective is to characterize the service rate policy that maximizes the infinite-horizon discounted reward and the long-run average reward. We show that for systems with an infinite buffer, the optimal service rate policy is monotone. However, the optimal policy is not necessarily monotone when the buffer size is finite. Furthermore, we prove that the set of possible optimal actions can be reduced to the lower boundary of the convex hull of the action space and develop an efficient policy iteration algorithm. Finally, we prove that the optimal service rate converges as the state tends to infinity. This property enables us to truncate the state space, facilitating the numerical computation of the optimal rate when the buffer space is infinite.
\end{abstract}
\setlength{\parindent}{2em}
\section{Introduction}
Customer abandonments are prevalent in many stochastic systems such as healthcare settings (see for example Batt and Terwiesh \cite{batt} and Armony et al. \cite{mor}), call centers (see for example Gong and Li \cite{gong} and Koole and Pot \cite{ger}), perishable inventory models (see for example Sarhangian et al. \cite{sarhangian} and Wu et al. \cite{wu}) etc. In spite of this wide range of applications, there are not many papers focusing on exact analysis and control of queues with abandonments. Our objective in this paper is to address this issue by considering service rate control in queues with abandonments.  

In particular, we consider a single stage Markovian queueing system with customer abandonments. We assume that customers arrive with respect to a Poisson process of rate $\lambda>0$ to a system whose buffer size could be infinite or finite (in which case customers finding the system full are lost). There is a single server whose service time is exponentially distributed. A service manager can choose the service mode of the server in a compact set $a\in \mathcal{A}$ or idle the server. When the server works in mode $a$, the service rate is $\mu(a)$ and  a  per unit time cost of $f(a)$ is incurred, where $\mu:\mathcal{A}\to\mathbb{R}^+$ and $f: \mathcal{A}\to\mathbb{R}$ are continuous functions in $\mathcal{A}$. The customers who are waiting in line are impatient and may leave after an exponential amount of time with rate $\theta>0$. When the server idles, all customers in the system are  subject to abandonments. There is a per unit time per customer holding cost $h\geq0$ incurred by customers in the system and a cost of $c\geq0$ per customer if a customer abandons the system. We consider two versions of the problem, in the first version, a reward $r\geq0$ is received at the time of customer arrivals whereas in the second version, reward $r$ is received at the time of service completions.  Our objective is to determine the structure of the optimal service rate control policy that maximizes the infinite-horizon discounted profit and the one that maximizes the long-run average profit for both versions of the problem.  To the best of our knowledge, this paper is the first to completely characterize the optimal service rate policy in $M/M/1 + M$ queues for both infinite and finite buffer systems under the infinite-horizon discounted and long-run average reward optimality criteria with a general continuous action space and service rate cost function. It also establishes that the set of optimal actions can be reduced to the lower boundary of the convex hull of the action space.

Service rate control has drawn attention in  queueing literature starting with Crabill \cite{crab1} who has shown that the optimal service rate is monotone in the number of customers in the system when it is chosen from a finite discrete set and there is a holding cost associated with the number of customers in the system. Crabill \cite{crab2} generalized this result to the setting where servers are subject to breakdowns. Using the techniques from their earlier paper (Weber and Stidham \cite{weber}), Stidham and Weber \cite{stidham} proved the optimality of monotone service rate policy when the holding cost is a convex increasing function of the number of customers in the system and the service rate can be chosen from a compact set. George and Harrison \cite{george} proved the same result with a more direct methodology that does not depend on dynamic programming. Jo and Stidham \cite{jo} showed that when the service time is a random sum of exponentially distributed phases, the optimal service rate is non-decreasing in the number of customers in the system for a given phase. Kumar et al. \cite{ravi} and Badian-Pessot et al. \cite{pamela} generalized the optimality of monotone service rate policy to single stage queues with Markov modulated arrivals and a removable server, respectively. 

There are several papers focusing on both admission control and service rate control in single stage Markovian queues. In an earlier paper, Serfozo \cite{dick} studied monotonic structures of birth-and-death processes, applied them to an M/M/s queue with controllable arrival rate and service rate, and identified sufficient conditions under which the optimal policies are monotonic. Ata and Sheneorson \cite{baris} proved that in a finite buffer system with a holding cost, a non-decreasing service cost, and reward associated with arrival, optimal arrival rate is non-increasing and the optimal service rate is non-decreasing in the number of customers in the system. In a similar setting, Adusumilli and Hasenbein \cite{john}
derived a simple, efficient algorithm for computing the optimal admission control and service rate  policy that minimizes the long-run average costs which include a holding cost, a service cost, and a rejection cost. Nelson and Kulkarni \cite{kulkarni} considered the same problem with multiple servers and proved that given a fixed service capacity, the optimal policy should use as many servers as possible with all busy servers serving at a common rate. Furthermore, they showed that the optimal arrival rate is again non-increasing and the optimal service rate is again non-decreasing in the number of customers in the system. More recently, Chen et al. \cite{chen} considered the problem of service rate control and (using an event-based approach) proved that the optimal policy is a bang-bang control. Finally, Zheng et al. \cite{gordon} studied the service rate control problem for an M/M/1 queue with server breakdowns in which the breakdown rate is assumed to be a function of the service rate and showed that the optimal policy is not necessarily monotone in the number of customers in the system.  

None of the papers discussed above considered optimal service rate control in the presence of abandonments. Down et al. \cite{Down} studied dynamic server control in a two-class service system with customer abandonments, demonstrating that the classical $c\mu$ rule fails to hold in general. Bhulai et al. \cite{sunjay} subsequently extended this problem to $K$ classes. In their seminal work, Atar et al. \cite{Atar} established the asymptotic optimality—in the many-server limiting regime—of the $c\mu/\theta$ rule for server scheduling in queues with abandonments. In a related line of research, Ko\c{g}a\u{g}a \cite{Kocaga} examined dynamic admission and service rate control for a $G/M/N+G$ queue via a diffusion model under the QED heavy-traffic regime. \c{S}im\c{s}ek et al. \cite{deniz}  proposed a deterministic approximation of the (mean) queueing dynamics (which they refer to  as a Unified Fluid Model) in queues where there is a dependency between the service times and abandonment times of customers. \c{C}ekyay \cite{bora} studied semi-Markov decision processes with unbounded transition rates and provided several sufficient conditions under which the value iteration procedure converges to the optimal value function and optimal deterministic stationary policies exist. In \cite{bora}, as an example of his methodology, \c{C}ekyay studied service rate control problem in queues with abandonments with the objective of minimizing the infinite-horizon discounted reward. However, his model requires certain conditions on the service cost function $f$, and the analysis under the long-run average reward criterion was left as a subject for future research. Using a more direct approach, we prove the optimality of monotone service rate policies in queues with infinite buffer space and abandonments both when customers pay at the time of arrival and at the time of service completion under the infinite horizon discounted optimality and the long-run average reward optimality criteria. Because the system is not uniformizable, we begin our analysis directly with the continuous-time Markov decision process formulation. After establishing the existence of optimal policies, we analyze the optimality equations analytically rather than relying on policy or value iteration to derive the structure of the optimal policy. Crucially, our approach applies to both optimality criteria and can be extended to other control problems, particularly those involving infinite state spaces and unbounded transition rates. Furthermore, we compare the optimal policies for two versions of the problem: payments at the time of arrival and payments at the time of service completion. Our results show that the optimal service rate policy is not necessarily monotone in finite buffer systems. Finally, we demonstrate that the set of possible optimal actions can be reduced to the lower boundary of the convex hull of the action space, and we develop a policy iteration algorithm to compute the optimal policy under both optimality criteria. 

The rest of the paper is organized as follows. Section \ref{formulation} formulates the problem for the infinite horizon discounted optimality criterion when customers pay at the time of arrival and Section \ref{S3} characterizes the optimal service rate policy for this version of the problem under both optimality criteria. Similarly, Section \ref{S4} focuses on the characterization of the optimal policies for the version where customers pay at the time of service completion. In Section \ref{S5}, we show that when the buffer size is finite, optimal service rate policy is not necessarily monotone. Using our structural results, in Section \ref{S6}, we develop an efficient policy iteration algorithm for finding the optimal service rate policy and use this algorithm in a numerical example. Section \ref{S7} concludes the paper. 
\section{Problem formulation} \label{formulation}
We start with the problem formulation for infinite buffer systems where customers pay at the time of arrival. Formulations for the other versions will be slight modifications of this original one. 

Consider the stochastic process $\{X^\pi(t):t\geq0\}$ under a policy $\pi\in\Pi$, where $\Pi$ denotes the set of all admissible service rate control policies and $X^\pi(t)=i\in S=\{0,1,2,\dots\}$ with $X^\pi(t)$ representing the total number of customers in the system at time $t$ under policy $\pi$. Let $A(i)$ be the set of possible actions in state $i\in S$. Then $A(0)=\{Idle\}$ and for $i\geq 1$, $A(i)=\mathcal{A}\cup\{Idle\}$ where $\mathcal{A}$ denotes the set of all possible working service modes available to the server. Let $r(i,a)$ be the reward earned per unit time  and $q_{ij}(a)$ be the transition rate of going from state $i$ to state $j$ for $i\neq j$ when action $a$ is chosen in state $i$. Then, the following happens in the time interval $[t,t+\mathrm{d}t]$:\par
(i) the system receives an infinitesimal reward $r(i,a)\mathrm{d}t$, and \par
(ii) a transition from state $i$ to state $j$ (with $i\neq j$) occurs with probability $q_{ij}(a)\mathrm{d}t+o(\mathrm{d}t)$; or the system remains at state $i$ with probability $1+q_{ii}(a)\mathrm{d}t+o(\mathrm{d}t)$,
where $q_{ii}(a):=-\sum_{j\neq i} q_{ij}(a)$.
Therefore, for  $i\geq 1$ and  $a\in \mathcal{A}$,

\[r(i,a)=\lambda r-hi-c(i-1)\theta-f(a),\]
and for $i\geq 0$ and $a=Idle$, 
\[r(i,Idle)=\lambda r-hi-c\theta i.\]
The rate matrix of the continuous time Markov chain when action $a\in\mathcal{A}$ is chosen in state $i\geq 1$ is

\[q_{ij}(a)=\left\{
\begin{array}{rcl}
    & \lambda & \mbox{for } j=i+1,\\
   &  \mu(a)+(i-1)\theta & \mbox{for } j=i-1,\\
  & -\lambda-\mu(a)-(i-1)\theta & \mbox{for } i=j,\\
  &0& \mbox{otherwise}.
\end{array}
\right.\]
and when action  $a=Idle$ is chosen in state $i\geq 0$, we have 
\[q_{ij}(Idle)=\left\{
\begin{array}{rcl}
    & \lambda & \mbox{for } j=i+1,\\
   &  i\theta & \mbox{for } j=\max\{0,i-1\},\\
  & -\lambda-i\theta & \mbox{for } i=j,\\
  &0& \mbox{otherwise}.
\end{array}
\right.\]
  To provide a uniform representation of the generator $q_{ij}(a)$ for $a\in A(i)$ for both actions in $\mathcal{A}$ and the idling action, we define $f: A\to \mathbb{R}$ with $f(Idle)=c\theta$ and $\mu: A\to \mathbb{R}^+$ with $\mu(Idle)=\theta$. Note that this assumption does not imply that idling provides service rate of $\theta$. Rather, when the server is idle, there is one additional customer who may abandon the system. Furthermore, under this assumption, the functions $\mu$ and $f$ remain continuous on $A$.
 \par
 
 The following result from Guo et al.\ \cite{survey} provides sufficient conditions for the continuous time Markov chain to be regular under any admissible policy.
 \begin{theorem}\label{A1}
For a continuous time MDP, if there exists a sequence $\{S_m:m\geq1\}$ of subsets of $S$, a nondecreasing function $\omega\geq1$ on $S$, constants $b_1\geq0$ and $b_2\neq 0$ such that \\
(i) $S_m\uparrow S $ and $\sup\{-q_{ii}(a):i\in S_m,a\in A\}<+\infty$ for each $m\geq 1$;\\
(ii) $\inf\{\omega(j):j\notin S_m\}\to+\infty$ as $m\to+\infty$;\\
(iii) $\sum_{j\in S}q_{ij}(a)\omega(j)\leq b_2 \omega(i)+b_1$ for all $(i,a)\in S\times A$.\\
Then the continuous time Markov chain under any admissible policy $\pi\in\Pi$ is regular.
 \end{theorem}
 In our model, since $\mu$ and $f$ are continuous functions and $A$ is compact, we know that  $\mu(A)$ is a compact set in $\mathbb{R}^+$ and  $f(A)$ is a compact set in $\mathbb{R}$. 
 Then by setting $S_m=\{0,1,\dots,m\}$, $\omega(j)=\max\{j,1\}$, $b_1=\lambda+\theta$ and $b_2=-\theta$, we can satisfy the three conditions of Theorem \ref{A1}.\par
 Next, for a policy $\pi\in\Pi$, let $\rho_t(B|i)$ be the probability that the policy $\pi$ chooses actions in a Borel set $B\in\mathcal{B}(A(i))$ at time $t$. Then, we define the reward function as
 \[R(t,i,\pi):=\int_{A(i)}r(i,a)\rho_t(\mathrm{d}a|i).\]
 \par
 
 Note that under our assumptions of customer abandonments and infinite buffer size, $q_{ii}(a)=-\lambda-\mu(a)-(i-1)\theta$ is not uniformly bounded for all states $i\in S$. As a result, the traditional uniformization method cannot be applied to transform the continuous time MDP into an equivalent discrete time MDP.  \c{C}ekyay \cite{bora} proposes a methodology that replaces the uniformization constant with a ``uniformization function'', thereby extending the uniformization method to CTMDPs with infinite state spaces. \cite{bora} shows that this methodology is applicable for the infinite-horizon discounted criterion. However, the resulting transition probabilities of the transformed discrete time Markov chain are more complicated, and the induction arguments for the value iteration becomes rather tedious.
 
 On the other hand, our approach works directly with the continuous time MDP and verifies the solution to the optimality equation through a limiting theorem. This provides a clearer idea of the structure of the optimal service rate control policy. 

 \par
 In the next section, we consider the infinite-horizon discounted reward and long-run average reward optimality in infinite buffer systems for the version where customers pay at the time of arrival and show that the optimal policy under both criteria is a monotone policy.

\section{Optimal policy when customers pay at arrival}\label{S3}
In Section 3.1, we characterize the structure of the optimal policy under the infinite horizon discounted optimality criterion whereas Section 3.2 is devoted to long-run average reward optimality.
\subsection{Infinite-horizon discounted reward optimality}\label{S31}
For the infinite-horizon discounted reward criterion with discount factor $\alpha>0$, if the system starts in state $i\in S$, then the value function under policy $\pi$ is defined as 
\[v(i,\pi,\alpha):=\mathbb{E}_i^\pi\big[\int_0^\infty e^{-\alpha t} R(t,X^\pi(t),\pi)\mathrm{d}t\big],\]
where the expectation is taken with respect to the continuous time Markov chain under policy $\pi$. The optimal discounted reward is then defined as
\[v^*(i,\alpha):=\sup_{\pi\in \Pi}v(i,\pi,\alpha).\]\par
In order to simplify the notation, whenever the discount factor $\alpha$ is not needed explicitly, we will write $v^*(i)$ instead of $v^*(i,\alpha)$ to denote the optimal infinite horizon discounted reward.\par 
Let $\mathbb{B}_\omega$ be the space of all real-valued functions on $S$ with 
\[\sup_{i\in S}\big\{\frac{|v(i)|}{\omega(i)}\big\}<+\infty.\]To guarantee that for all 
$i\in S$, $v^*(i)<\infty$ and there exists a unique solution to the optimality equation that belongs to $\mathbb{B}_\omega$, we need to verify the following conditions of Guo et al. \cite{survey}.
\begin{theorem}\label{A2}
For the infinite-horizon discounted reward criterion, if\\
    (i) for every $(i,a)\in S\times A$ and some constant $M>0$,
\[|r(i,a)|\leq M\omega(i),\]
where $\omega(i)$ is as defined in Theorem \ref{A1},\\
(ii) the discount factor $\alpha>0$ satisfies $\alpha>b_2$ where $b_2$ is as defined in Theorem \ref{A1},\\
then the value function $v^*(i)=\sup_{\pi\in \Pi}v(i,\pi)\in \mathbb{B}_\omega$ and satisfies
\begin{equation}\label{OE2}
   \alpha v^*(i)=\sup_{a\in A(i)}\{r(i,a)+\sum_{j\in S}q_{ij}(a)v^*(j)\}\qquad \forall\  i\in S.
\end{equation}
\end{theorem}
Since $f(A)$ is a compact set in $\mathbb{R}$, we know that there exists $f_+$ such that $|f(a)|\leq f_+$ for all $a\in A$. Then, we set $M=\lambda +r+f_+ +h+c\theta$ and since $b_2=-\theta<0<\alpha$, the conditions in Theorem \ref{A2} are satisfied and the optimal value function satisfies the optimality equations in (\ref{OE2}). \par
Let $\mathbb{F}$ be the set of all stationary deterministic policies. We will show that we only need to search for the optimal policy in $\mathbb{F}$ because the following conditions of Guo et al. \cite{survey} hold for our model.
\begin{theorem}\label{A3}
For the continuous time MDP, if\\ 
    (i) the action set $A(i)$ is compact for each $i\in S$,\\
(ii) the function $r(i,a),q_{ij}(a)$ and $\sum_{k\in S}q_{ik}(a)\omega(k)$ are all continuous in $a\in A(i)$ for each fixed $i,j\in S$,\\
(iii) there exists a nonnegative function $\omega^\prime$ on $S$, and constants $b_2^\prime >0,\ b_1^\prime\geq0$ and $M^\prime>0$ such that for all $i\in S$
\[q(i)\omega(i)\leq M^\prime \omega^\prime(i)\]
where $q(i)=\sup_{a\in A(i)}\{-q_{ii}(a)\}$, and 
\[\sum_{j\in S}\omega^\prime(j)q_{ij}(a)\leq b_2^\prime \omega^\prime(i)+b_1^\prime \mbox{ for all $(i,a)\in S\times A$,}\]
then there exists a stationary deterministic policy that attains the maximum infinite-horizon discounted reward.
\end{theorem}
Note that (i) and (ii) follow immediately from our assumptions.  Define $\mu_-=\min\{\mu(A)\}$ and $\mu_+=\max\{\mu(A)\}$. Then for (iii), we can choose $\omega^\prime(i)=\max\{i^2,1\}$, $M^\prime=\lambda+\mu_++\theta$, $b_1^\prime=3\lambda$ and $b_2^\prime=\lambda$. Thus, (i), (ii) and (iii) are all satisfied and there exists a stationary deterministic policy that attains the maximum infinite-horizon discounted reward. \par
For a stationary deterministic policy $\pi\in\mathbb{F}$, if the system starts in state $i$ and the discount factor is $\alpha$, then the value function can be simplified to
\[v(i,\pi,\alpha)=\mathbb{E}_i^\pi\big[\int_0^\infty e^{-\alpha t} r(X^\pi(t),a_{X^\pi(t)})\mathrm{d}t\big],\]
where $a_{X^\pi(t)}$ denotes the action taken in state $X^\pi(t)$ at time $t$ under policy $\pi$ and the expectation is taken with respect to the continuous time Markov chain under policy $\pi$.

We begin by stating the main result of this subsection. Specifically, we prove that the optimal service rate control policy is monotone under the infinite-horizon discounted reward criterion.
\begin{theorem}\label{T33}
    The optimal policy $\pi^*$ that maximizes the infinite-horizon discounted reward is monotone non-decreasing in the number of customers in the system.
\end{theorem}

The remainder of this subsection is dedicated to proving this main result.
For an optimal policy $\pi^*\in\mathbb{F}$, let $a_i^*$ be the optimal action in state $i$. From the optimality equation in (\ref{OE2}), we know that the optimal policy satisfies the equations
\begin{eqnarray}\label{OE1}
\alpha v^*(0)=\lambda r-\lambda v^*(0)+\lambda v^*(1),\qquad \nonumber\\
       \alpha v^*(i)=r(i,a_i^*)+\sum_{j\in S}q_{ij}(a_i^*)v^*(j)\qquad \forall\  i\geq 1.
\end{eqnarray}
Then, we can expand the equations in (\ref{OE1}) as follows

\small
\begin{equation}\label{EEv1}
    \alpha v^*(i)=\left\{
\begin{array}{rcl}
    & \lambda r+\lambda v^*(1)-\lambda v^*(0) & \mbox{for }i=0,\\
    & \lambda r-hi-c(i-1)\theta-f(a_i^*)+\lambda v^*(i+1)\\
    &+(\mu(a_i^*)+(i-1)\theta)v^*(i-1)-(\lambda+\mu(a_i^*)+(i-1)\theta)v^*(i)&\mbox{for }i>0.\\
\end{array}
\right.
\end{equation}\par
\normalsize
For any vector $v$ and $i\geq 0$, define $\Delta v(i)=v(i+1)-v(i)$. Then from (\ref{EEv1}) with some algebra we have
\normalsize
\begin{equation}\label{OE311}
     \alpha\Delta v^*(i)=\left\{
\begin{array}{rcl}
     &-h-f(a_1^*)+\lambda\Delta  v^*(1)-(\lambda+\mu(a_1^*))\Delta  v^*(0)& \mbox{for } i=0,  \\
     & -h-c\theta+\lambda\Delta  v^*(i+1)-(\lambda+i\theta)\Delta  v^*(i)+(i-1)\theta\Delta  v^*(i-1)\\
     &-f(a_{i+1}^*)+f(a_i^*)-\mu(a_{i+1}^*)\Delta  v^*(i)+\mu(a_i^*)\Delta v^*(i-1) & \mbox{for } i>0.\\
     
\end{array}
\right.
\end{equation}
\normalsize

Similarly, define, for $i\geq 0$, $\Delta^2 v(i)=\Delta v(i+1)-\Delta v(i)$ and for $i\geq 1$, $w(i)=-f(a_{i+1}^*)+f(a_i^*)-\mu(a_{i+1}^*)\Delta  v^*(i)+\mu(a_i^*)\Delta v^*(i-1)$, we then have
\begin{equation}\label{OE31}
\alpha\Delta^2 v^*(i)=\left\{
\begin{array}{rcl}
    & f(a_1^*)-c\theta+\lambda\Delta^2 v^*(1)-(\lambda+\theta)\Delta^2 v^*(0)  \\
    &+(\mu(a_1^*)-\theta)\Delta v^*(0)+w(1)&\mbox{for } i=0,\\
    & \lambda\Delta^2 v^*(i+1)-(\lambda+(i+1)\theta)\Delta^2 v^*(i)+(i-1)\theta\Delta^2 v^*(i-1)\\
    & +w(i+1)-w(i)&\mbox{ for } i>0.
\end{array}
\right.
\end{equation}
\par
Before presenting the formal proof of our main result, we outline the key steps and technical arguments to provide a high-level overview of our approach.
\begin{enumerate}
    \item We show that $\Delta v^*(i-1)$, together with the $\mu$ and $f$ functions, completely characterizes the optimal service rate control policy in state $i$ (Lemma \ref{Lemma31}). Consequently, $\mu(a_i^*)$ is non-decreasing in $i$ if and only if $\Delta^2 v^*(i)$ is less than or equal to $0$ (Corollaries \ref{C33} and  \ref{C33.5}). 
    \item We show that for $i\geq1$, $w(i)=-\gamma_{i-1}\Delta^2 v^*(i-1)$ for some constant vector $\{\gamma_i\}_{i=0}^\infty$ (Lemmas \ref{L32} and \ref{C35}). As a result, we can rewrite equation (\ref{OE31}) as a system of equations for $\Delta^2 v^*(i),\ i\geq0$.
    \item We truncate this system of equations and consider the first $L$ number of equations (equation (\ref{EqL})). We show that the solution to this truncated system of equations $\Delta^2 v_L(i), 0\leq i\leq L$ is less than or equal to $0$ (Lemmas \ref{L37} and  \ref{L36}).
    \item We show that $\Delta^2 v_L(i)$ is monotone non-increasing in $L$ and pointwise converges to $\Delta^2 v^*(i)$. As a result, the solution to the original system of equations is less than or equal to $0$ (Corollary \ref{C34}).
    \item We conclude that the optimal service rate is monotone non-decreasing in the number of customers in the system (Corollary \ref{CT33}).
\end{enumerate}

First, we prove the following lemma which characterizes the set of optimal actions in each state.
\begin{lemma}\label{Lemma31}
Under infinite horizon discounted optimality criterion, for $i\geq 1$, $a_i^*$ satisfies
\[a_i^*\in\mathop{argmax}_{a\in A}\{-f(a)-\mu(a)\Delta v^*(i-1)\}.\]
\end{lemma}
\begin{proof}
    From the optimality equation in (\ref{OE2}), with some algebra, for $i\geq 1$,  we have
    \begin{eqnarray}
        \alpha v^*(i)&=&\sup_{a\in A}\{\lambda r-ih-c(i-1)\theta-f(a)+\lambda \Delta v^*(i)-((i-1)\theta+\mu(a))\Delta v^*(i-1)\}\nonumber\\
        &=&\lambda r-ih-c(i-1)\theta+\lambda \Delta v^*(i)-(i-1)\theta \Delta v^*(i-1)+\sup_{a\in A}\{-f(a)-\mu(a)\Delta v^*(i-1)\}.\nonumber
    \end{eqnarray}
    Since $a$ is chosen in a compact set $A$ and $f(a)$ is continuous, the supremum is attained and we have, for $i\geq 1$,
    \[a_i^*\in\mathop{argmax}_{a\in A}\{-f(a)-\mu(a)\Delta v^*(i-1)\}.\]
\end{proof}
We then immediately have the following characterization of $a_1^*$.
\begin{corollary}\label{C31}
Optimal action in state $1$ satisfies,
    \[f(a_1^*)-c\theta+(\mu(a_1^*)-\theta)\Delta v^*(0)\leq0.\]
\end{corollary}
\begin{proof}
    From Lemma \ref{Lemma31}, we have
    \[a_1^*\in\mathop{argmax}_{a\in A}\{-f(a)-\mu(a)\Delta v^*(0)\}.\]
    Thus, for $a=Idle$, 
    \[f(a_1^*)+\mu(a_1^*)\Delta v^*(0)\leq c\theta+\theta\Delta v^*(0)\]
    and the result follows.
\end{proof}
The following bounds on $w(i)$ will be used to express it in terms of $\Delta^2v^*(i-1)$.
\begin{lemma}\label{L32}
    For $i\geq1 $, we have 
    \begin{equation}\label{Core}
        -\mu(a_i^*)\Delta^2v^*(i-1)\leq w(i)\leq -\mu(a_{i+1}^*)\Delta^2 v^*(i-1).
    \end{equation}
\end{lemma}
\begin{proof}
    From Lemma \ref{Lemma31}, we have 
    \[-f(a_{i+1}^*)-\mu(a^*_{i+1})\Delta v^*(i)=\max_{a\in A}\{-f(a)-\mu(a)\Delta v^*(i)\}\geq -f(a^*_i)-\mu(a^*_{i})\Delta v^*(i).\]
    Then,
    \begin{eqnarray}\label{Iq7}
        w(i)&=&-f(a^*_{i+1})-\mu(a^*_{i+1})\Delta v^*(i)+f(a^*_i)+\mu(a^*_i)\Delta v^*(i-1) \nonumber\\
        &\geq& -f(a^*_i)-\mu(a^*_{i})\Delta v^*(i)+f(a^*_i)+\mu(a^*_i)\Delta v^*(i-1) \nonumber\\
        &=&-\mu(a^*_i)\Delta^2 v^*(i-1).
    \end{eqnarray} 
    Similarly, from Lemma \ref{Lemma31}, we have
    \[-f(a^*_i)-\mu(a^*_i)\Delta v^*(i-1)=\max_{a\in A}\{-f(a)-\mu(a)\Delta v^*(i-1)\}\geq -f(a^*_{i+1})-\mu(a^*_{i+1})\Delta v^*(i-1).\]
    Then,
    \begin{eqnarray}\label{Iq8}
        w(i)&=&-f(a^*_{i+1})-\mu(a^*_{i+1})\Delta v^*(i)+f(a^*_i)+\mu(a^*_i)\Delta v^*(i-1)\nonumber\\
        &\leq& -f(a^*_{i+1})-\mu(a^*_{i+1})\Delta v^*(i)+f(a^*_{i+1})+\mu(a^*_{i+1})\Delta v^*(i-1)\nonumber\\
        &=&-\mu(a^*_{i+1})\Delta^2 v^*(i-1). 
    \end{eqnarray}
    Inequalities (\ref{Iq7}) and (\ref{Iq8}) together yield (\ref{Core}).
\end{proof}
The following corollaries follow immediately from Lemma \ref{L32} and provide properties of $\mu(a_i^*)$ under certain conditions.
\begin{corollary}
    If $\Delta^2 v^*(i-1)<0$ for some $i\geq 1$, then $\mu(a_i^*)\leq \mu(a_{i+1}^*)$.
\end{corollary}
\begin{corollary}\label{C33}
     If $\Delta^2 v^*(i-1)<0$ for all $i\geq1$, then $\mu(a_i^*)$ is monotone non-decreasing in $i$, for $i\geq 1$. Thus, the optimal policy has a monotone structure.
\end{corollary}
\begin{corollary}\label{C33.5}
    If $\Delta^2 v^*(i-1)=0$ for some $i\geq 1$, then there exists an optimal policy $\pi^*$ under which  $\mu(a_i^*)=\mu(a_{i+1}^*)$.
\end{corollary}\par

The following lemma follows immediately from equation (\ref{Core}) and the Intermediate Value Theorem (see Theorem 4.23 in \cite{rudin}).
\begin{lemma}\label{C35}
    For $i\geq 1$, there exists some $\gamma_{i-1}\in[\min\{\mu(a_{i}^*),\mu(a_{i+1}^*)\},\max\{\mu(a_{i}^*),\mu(a_{i+1}^*)\}]$, such that $w(i)=-\gamma_{i-1}\Delta^2 v^*(i-1)$.
\end{lemma}

\par
Furthermore, in order to simplify the notation, we set \[G=\frac{f(a_1^*)-c\theta+(\mu(a_1^*)-\theta)\Delta v^*(0)}{\lambda}\leq0\]
where the inequality follows from Corollary $\ref{C31}$.\par Plugging these into the equations in (\ref{OE31}),  we have
\begin{equation}\label{OE40}
    (\alpha+\lambda+\theta+\gamma_0)\Delta^2 v^*(0)=
     \lambda\Delta^2 v^*(1)+\lambda G,
\end{equation}
\begin{equation}\label{OE41}
    (\alpha+\lambda+(i+1)\theta+\gamma_i)\Delta^2 v^*(i)=\lambda\Delta^2 v^*(i+1)+\left(\left(i-1\right)\theta+\gamma_{i-1}\right)\Delta^2 v^*(i-1)\mbox{ for } i>0.
\end{equation}
\par
In order to simplify the notation and derive a clear relationship between $\Delta^2 v^*(i)$ and $\Delta^2 v^*(0)$, for $i\geq0$, define $p_i=\frac{\alpha+(i+1)\theta+\gamma_i}{\lambda}$. Similarly, let $q_0=p_0+1$, $q_1=p_1q_0+q_0-p_0+\frac
{\alpha+\theta}{\lambda}$ and define $q_{i+2}=q_{i+1}p_{i+2}+q_{i+1}-p_{i+1}q_i+q_i\frac{\alpha+\theta}{\lambda}$ for $i\geq 0$.\par
In the remainder of this subsection, we show that the solution to the system of equations \eqref{OE40}–\eqref{OE41} is non-positive.

The following lemma states the relationship between $\Delta^2 v^*(i)$ for $i\geq1$ and $\Delta^2 v^*(0)$.
\begin{lemma}\label{L33}
    For $i\geq0$,
    \[\Delta^2 v^*(i+1)=q_i\Delta^2 v^*(0)-\Pi_{j=1}^{i}(p_j+1)G\]
    with the convention that $\Pi_{j=1}^0(p_j+1)=1$.
    
\end{lemma}
\begin{proof}
Note that for $i\geq0$, $\alpha+\lambda+(i+1)\theta+\gamma_i=\lambda(p_i+1)$ and $i\theta+\gamma_i=\lambda(p_i-\frac{\alpha+\theta}{\lambda}).$ For $i=0$, from equation (\ref{OE40}), we have 
\[\lambda(p_0+1)\Delta^2 v^*(0)=\lambda\Delta^2 v^*(1)+\lambda G,\]
and hence,
\[\Delta^2 v^*(1)=(p_0+1)\Delta^2 v^*(0)-G=q_0\Delta^2 v^*(0)-G.\]
For $i=1$, from equation (\ref{OE41}), we have
\[\lambda(p_1+1)[q_0\Delta^2 v^*(0)-G]=\lambda\Delta^2 v^*(2)+\lambda(p_0-\frac{\alpha+\theta}{\lambda})\Delta^2 v^*(0),\]
and hence,
\[\Delta^2 v^*(2)=(p_1q_0+q_0-p_0+\frac
{\alpha+\theta}{\lambda})\Delta^2 v^*(0)-(p_1+1)G=q_1\Delta^2 v^*(0)-(p_1+1)G.\]
We use induction to prove the result for $i>1$. Suppose that the result holds for all $k\leq i$. For $k=i+1$, from equation (\ref{OE41}) we have 
\[\lambda(p_{i+1}+1)[q_i\Delta^2 v^*(0)-\Pi_{j=1}^{i}(p_j+1)G]=\lambda\Delta^2 v^*(i+2)+\lambda(p_{i}-\frac{\alpha+\theta}{\lambda})q_{i-1}\Delta^2 v^*(0),\]
and thus,
\begin{eqnarray}
    \Delta^2 v^*(i+2)&=&(p_{i+1}q_i+q_i-p_{i}q_{i-1}+q_{i-1}\frac
{\alpha+\theta}{\lambda})\Delta^2 v^*(0)-\Pi_{j=1}^{i+1}(p_j+1)G\nonumber\\
&=&q_{i+1}\Delta^2 v^*(0)-\Pi_{j=1}^{i+1}(p_j+1)G,\nonumber
\end{eqnarray}
which concludes the proof.
\end{proof}
The next lemma summarizes the relationship between $q_i,p_i,q_{i+1}$ and $p_{i+1}$ for $i\geq0$.
\begin{lemma}\label{L34}
For $i\geq0$,
\[0<q_ip_{i+1}<q_{i+1}<q_i(p_{i+1}+1).\]
\end{lemma}
\begin{proof}
Note that $p_i>\frac{\alpha+\theta}{\lambda}>0$ for all $i\geq 0$. We have $q_0>0$, 
\[q_1=(p_1q_0+q_0)-p_0+\frac{\alpha+\theta}{\lambda}<q_0(p_1+1), \mbox{ and}\]
\[q_1=p_1q_0+(q_0-p_0)+\frac{\alpha+\theta}{\lambda}=p_1q_0+1+\frac{\alpha+\theta}{\lambda}>q_0p_1>0.\]
Thus, the result holds for $i=0$. 
We use induction to prove the result for $i>0$. Assume the result holds for $k=0,\dots,i-1$. We show that it also holds for $k=i$. From the induction hypothesis we have $q_{i-1}>0$ and $q_i>p_iq_{i-1}$. Then 
\[q_{i+1}=q_{i}p_{i+1}+q_{i}-(p_{i}-\frac{\alpha+\theta}{\lambda})q_{i-1}<q_{i}(p_{i+1}+1)\]
and \[q_{i+1}=q_{i}p_{i+1}+(q_{i}-p_{i}q_{i-1})+\frac{\alpha+\theta}{\lambda}q_{i-1}>q_{i}p_{i+1}>0.\]
Therefore, the result holds.
\end{proof}
We next focus on the computation of $\Delta v^*(i),$ for all $i\geq0$ (as given in equations (\ref{OE40}) and (\ref{OE41})). In order to do this, we start with the following set of equations for a truncated system. For $L\geq1$, define
\small
\begin{equation}\label{EqL}
    (\alpha+\lambda+(i+1)\theta+\gamma_i)\Delta^2 v_L(i)=\left\{
\begin{array}{rcl}
     &\lambda G+\lambda\Delta^2 v_L(1)& \mbox{for } i=0,  \\
     & \lambda\Delta^2 v_L(i+1)+((i-1)\theta+\gamma_{i-1})\Delta^2 v_L(i-1) & \mbox{for } 0<i\leq L,\\
     &0&\mbox{for } i>L.
     
\end{array}
\right.
\end{equation}
\normalsize
We immediately have the following result from Lemma \ref{L33}.
\begin{lemma}\label{L35}
   For $0\leq i\leq L$, 
\[\Delta^2 v_L(i+1)=q_i\Delta^2 v_L(0)-\Pi_{j=1}^{i}(p_j+1)G,\]
    with the convention that $\Pi_{j=1}^0(p_j+1)=1$.
\end{lemma}
Note that in (\ref{EqL}), we have $\Delta^2 v_L(L+1)=0$. Then from Lemma \ref{L35}, 
\[q_{L}\Delta^2 v_L(0)-\Pi_{j=1}^L(p_j+1)G=0\]
and
\begin{equation}\label{EqL0}
    \Delta^2 v_L(0)=\Pi_{j=1}^L(p_j+1)\frac{G}{ q_L}\leq0
\end{equation}
where the inequality follows from Lemma \ref{L34} and the fact that $G\leq0$. Then for $1<i\leq L$,
\begin{equation}\label{EqL1}
    \Delta^2 v_L(i)=\frac{q_{i-1}}{q_L}\Pi_{j=1}^L(p_j+1)G-\Pi_{j=1}^{i-1}(p_j+1)G.
\end{equation}
\par
Since, $q_L>0$, we know that (\ref{EqL0}) and (\ref{EqL1}) constitute a unique solution to the system of equations in (\ref{EqL}). We then have the following results.
\begin{lemma}\label{L37}
     Let $\Delta^2 v_L(i),i\in S$ be the unique solution to the system of equations in (\ref{EqL}). If $G=0$, then $\Delta^2 v_L(i)=0$ for all $i\in S$ and $L\geq 1$.
\end{lemma}
\begin{lemma}\label{L36}
    Let $\Delta^2 v_L(i),i\in S$ be the unique solution to the system of equations in (\ref{EqL}). If $G<0$, then\\
    (i) $\Delta^2 v_L(i)<0$ for all $i\leq L$ and\\
    (ii) $\Delta^2v_L(i)>\Delta^2v_{L+1}(i)$ for all $i\leq L$.
    
\end{lemma}
\begin{proof}
    (i) Since $G<0$, (\ref{EqL0}) shows that $\Delta^2 v_L(0)<0$. Furthermore, since $q_{i-1}(p_{i}+1)>q_i$ (see Lemma \ref{L34}), we have 
    \[\Pi_{j=i}^L(p_j+1)q_{i-1}>\Pi_{j=i+1}^L(p_j+1)(p_i+1)q_{i-1}>\Pi_{j=i+1}^L(p_j+1)q_{i}>\cdots>(p_L+1)q_{L-1}>q_L.\] Therefore,
    \[\Delta^2 v_L(i)=\Pi_{j=1}^{i-1}(p_j+1)G[\Pi_{j=i}^L(p_{j}+1)\frac{q_{i-1}}{q_L}-1]<0.\]
    (ii) For $i\leq L$, we have 
    \begin{eqnarray}
       \Delta^2v_L(i)&=& \frac{q_{i-1}}{q_L}\Pi_{j=1}^L(p_j+1)G-\Pi_{j=1}^{i-1}(p_j+1)G\nonumber\\
       &>&\frac{q_{i-1}}{q_L}\frac{q_L(p_{L+1}+1)}{q_{L+1}}\Pi_{j=1}^L(p_j+1)G-\Pi_{j=1}^{i-1}(p_j+1)G\nonumber\\
       &=&\frac{q_{i-1}}{q_{L+1}}\Pi_{j=1}^{L+1}(p_j+1)G-\Pi_{j=1}^{i-1}(p_j+1)G\nonumber\\
       &=&\Delta^2 v_{L+1}(i),\nonumber
    \end{eqnarray}
    where the inequality follows from Lemma \ref{L34} and the assumption that $G<0$.
\end{proof}
\begin{corollary}\label{C34}
    Let $\Delta^2 v_L(i),i\in S$ be the unique solution to the system of equations in (\ref{EqL}), and $v^*(i)$ be the optimal infinite-horizon expected discounted reward. Then for each state $i\in S$
    \[\lim_{L\to+\infty}\Delta^2 v_L(i)=\Delta^2 v^*(i).\]
\end{corollary}
\begin{proof}
     If $G=0$, the result follows immediately from Lemma \ref{L37}. Assume $G<0$. From Theorem \ref{A2}, we know that the unique solution to the optimality equations in (\ref{OE1}) belongs to $\mathbb{B}_\omega$. Fix $i\in S$, from Lemma \ref{L36} (ii) we know that when $L\geq i$, $\Delta^2 v_L(i)$ satisfying the truncated optimality equations is monotone non-increasing with respect to $L$. Therefore, it pointwise converges to $\Delta^2 v^*(i)$. 
\end{proof}

\begin{corollary}\label{CT33}
    For all $i\in S$, $\Delta^2 v^*(i)\leq 0$. 
\end{corollary}

\begin{proof}
    From Lemmas \ref{L37} and \ref{L36} and Corollary \ref{C34}, we know that if $G<0$, then $\Delta^2 v^*(i)<0$  and if $G=0$, then $\Delta^2 v^*(i)=0$. 
\end{proof}
Hence, from Corollary \ref{C33},  if $G<0$, the optimal policy is monotone non-decreasing with respect to $i\in S$. Similarly, if $G=0$, from Corollary \ref{C33.5}, we know that there exists an optimal policy that idles the server in all states. This concludes the proof of Theorem \ref{T33}.

\begin{remark}\label{R31}
Note that the structure of the optimal service rate policy in Theorem \ref{T33} matches that of classical systems without abandonments (see, for example, George and Harrison \cite{george}, Kumar et al. \cite{ravi}, and  Stidham and Weber \cite{stidham}). The intuition behind this similarity is that the service rate control policy cannot influence the arrival process. That is, the policy cannot directly dictate customer abandonment behavior; it can only control the departure rate (via service completions) of customers already present in the queueing system. Consequently, an optimal policy still aims to minimize system congestion, just as it would in a setting without abandonments. However, a key distinction of our work is that we do not require $f$ to be non-decreasing (as in \cite{george}, \cite{ravi}, and \cite{stidham}) or convex (as in \cite{ravi}) to establish this monotonicity result. Since our proof readily extends to the case where $\theta = 0$, it also provides an alternative way to establish the optimality of monotone service rate policies in infinite-buffer systems without needing these additional assumptions on the cost function $f$.
\end{remark}

\subsection{Properties of the infinite horizon discounted reward optimal policy}\label{S32}
In this section, we  study properties of the optimal expected discounted reward $v^*(i)$ and the optimal policy $\pi^*$.\par
Since for any state $i$ and action $a$, $r(i,a)$ is the sum of the reward $\lambda r$ and the cost $C(i,a)=-hi-c(i-1)\theta-f(a)$, for $i\geq 1$, and $C(0,Idle)=0$, we can split the value function into two parts. For a stationary deterministic policy $\pi\in\mathbb{F}$, if the system starts in state $i$ with discount factor $\alpha$, the value function can be written as
\[v(i,\pi,\alpha)=\mathbb{E}_{i}^\pi\big[\int_0^\infty e^{-\alpha t}(\lambda r+C(X^\pi(t),a_{X^\pi(t)}))\mathrm{d}t\big]=\frac{1}{\alpha}\lambda r+\mathbb{E}_{i}^\pi\big[\int_0^\infty e^{-\alpha t}C(X^\pi(t),a_{X^\pi(t)})\mathrm{d}t\big].\]
Therefore, $v^*(i)=\sup_{\pi\in\Pi}\{v(i,\pi,\alpha)\}=\frac{1}{\alpha}\lambda r+\sup_{\pi\in\Pi}\mathbb{E}_{i}^\pi\big[\int_0^\infty e^{-\alpha t}C(X^\pi(t),a_{X^\pi(t)})\mathrm{d}t\big]$ and we immediately have the following theorem.
\begin{Proposition}\label{Pro31}
    The infinite horizon discounted reward optimal policy $\pi^*$ does not  depend on $r$.
\end{Proposition}
Proposition \ref{Pro31} implies that the total expected reward received from the arriving customers is the same under all service rate control policies. We next prove that the optimal actions at any state $i\geq1$ can only be chosen from the lower boundary of the convex hull of the action space.

From Lemma \ref{Lemma31} we know that, for  $i\geq 1$,
\[a_i^*\in\mathop{argmax}_{a\in A}\{-f(a)-\mu(a)\Delta v^*(i-1)\}.\]
Therefore, consider the following optimization problem
\begin{equation}\label{NL}
    \max_{(x,y)\in\{(\mu(a),f(a)), a\in A\}}\qquad -y-x\Delta v^*(i-1).
\end{equation}
Since $A$ is a compact set and functions $f$ and $\mu$ are continuous, we know that $H=\{(\mu(a),f(a)),a\in A\}$ is a compact set in $\mathbb{R}^2$. Note that the objective function of the optimization problem is linear. As a result, we have the following theorem.
\begin{theorem}\label{T34CH}
     All optimal actions $a_i^*$ in state $i\geq1$ are on the lower boundary of $H=\{(\mu(a),f(a)),a\in A\}$. Furthermore, there exists an optimal policy where the optimal actions $a_i^*$ are the extreme points of the convex hull of $H$ in $\mathbb{R}^2$. That is to say, let $conv(H)$ be the convex hull of $H$, we have for all $i\geq 1$,\\
     (i) $f(a_i^*)=\inf\{f(a):a\in A,\mu(a)=\mu(a_i^*)\}$ and\\
    (ii) there exists an $a_i^*$ such that $(\mu(a_i^*),f(a_i^*))$ is an extreme point of $conv(H)$.
    
\end{theorem}
\begin{proof}
     Since the objective function of the optimization problem (\ref{NL}) is linear and $H$ is compact (thus a closed bounded set in $\mathbb{R}^2$), we know that there always exists a solution $(x^*,y^*)$ to (\ref{NL}) which is an extreme point of the $conv(H)$ (see Corollary 32.3.2 in Rockafellar \cite{convex}). Furthermore, in the objective function, the coefficient of $y$ ($f(a)$) is $-1$. As a result, all solutions to (\ref{NL}) must be on the lower boundary of $conv(H)$.
\end{proof}
Theorem \ref{T34CH} describes a set that contains all possible actions for an optimal policy. In general, it is  not easy to find every extreme point of $conv(H)$. However, in practice, we only need to consider the lower boundary of $conv(H)$. Note that this lower boundary can be written as a function from $\mu(A)$ to $f(A)$. Furthermore, let $f^*$ be the lower boundary of $conv(H)$, then $f^*$ is continuous, semi-differentiable, and convex.
The domain of $f^*$ is a closed interval $[\mu_-,\mu_+]$ and it contains all possible optimal actions of the service rate control problem. Although analytically characterizing the lower boundary of the convex hull $f^*$ may be difficult in practice, the main contribution of this result is that it demonstrates that the standard assumptions of convexity and monotonicity on the service cost function $f$ are not required to obtain monotone optimal policies (see Remark \ref{R31}). Furthermore, this structural insight provides a powerful filtering mechanism, effectively eliminating service actions that can never be optimal. For example,
if one uses the policy iteration algorithm (see Section 6.4 of Puterman \cite{put}) to find the optimal policy, instead of searching inside the whole action space $A$ and solving (\ref{NL}), we can use convex optimization methods, for instance, gradient descent, to solve (\ref{NL}) on $f^*$, which is guaranteed to find the optimal action and accelerate the speed of each iteration of the policy iteration algorithm. In other words, the optimization problem in (\ref{NL}) is transformed into the convex optimization problem
\begin{equation}\label{NL1}
    \max_{\mu\in[\mu_-,\mu_+]}\qquad -f^*(\mu)-\mu\Delta v^*(i-1).
\end{equation}


Next, we show that the set of possible optimal actions can be narrowed down even further. First, we obtain bounds on $v^*$ and $\Delta v^*$.

\begin{Proposition}\label{P31}
    For the infinite-horizon discounted reward criterion, for all  $i\in S$, we have
\[v^*(i)\geq \frac{\lambda}{\alpha}(r-\frac{h+c\theta}{\alpha+\theta})-\frac{h+c\theta}{\alpha+\theta}i\]    and \[ \Delta v^*(i)\geq -\frac{h+c\theta}{\alpha+\theta}.\]
Furthermore, if  $f\geq0$, for all  $i\in S$, we have
    \[\frac{\lambda}{\alpha}(r-\frac{h+c\theta}{\alpha+\theta})-\frac{h+c\theta}{\alpha+\theta}i\leq v^*(i)<\frac{\lambda}{\alpha}r\]
 and    \[-\frac{h+c\theta}{\alpha+\theta}\leq \Delta v^*(i)<0.\]
\end{Proposition}
\begin{proof}
    Consider the policy $\pi^\prime$ that idles the server in all states. We have 
\begin{equation}\label{Bd1}
    v(i,\pi^\prime,\alpha)=\frac{\lambda}{\alpha}(r-\frac{h+c\theta}{\alpha+\theta})-\frac{h+c\theta}{\alpha+\theta}i.
\end{equation}
Therefore, 
\[v^*(i)\geq \frac{\lambda}{\alpha}(r-\frac{h+c\theta}{\alpha+\theta})-\frac{h+c\theta}{\alpha+\theta}i.\]
Since $\Delta^2 v^*(i)\leq0$, $\Delta v^*(i)$ is monotone non-increasing in $i$. If for some $i\in S$, $\Delta v^*(i)<-\frac{h+c\theta}{\alpha+\theta}$, then
$\Delta v^*(j)<-\frac{h+c\theta}{\alpha+\theta}$ for all $j\geq i$. Therefore, when $j$ is large enough, we have $v^*(j)<\frac{\lambda}{\alpha}(r-\frac{h+c\theta}{\alpha+\theta})-\frac{h+c\theta}{\alpha+\theta}j$, which contradicts  (\ref{Bd1}). Thus, we must have  $\Delta v^*(i)\geq -\frac{h+c\theta}{\alpha+\theta}$ for all $i\geq 0$.

Furthermore, if  $f\geq 0$, then the cost incurred in state $i\geq 1$ under action $a$, $C(i,a)<0$ for all $i\geq 1$. Therefore, for all $i\geq 0$, we immediately have $v(i,\pi,\alpha)<\frac{1}{\alpha}\lambda r$ and
\[\frac{\lambda}{\alpha}(r-\frac{h+c\theta}{\alpha+\theta})-\frac{h+c\theta}{\alpha+\theta}i\leq v^*(i)<\frac{\lambda}{\alpha}r.\]
Similarly, plugging this inequality into equation ($\ref{EEv1}$) for $i=0$, we have
\[-\frac{h+c\theta}{\alpha+\theta}\leq\Delta v^*(0)<0.\]  
    Since $\Delta^2 v^*(i)\leq0$,  $\Delta v^*(i)\leq \Delta v^*(0)<0$. Then, for all $i\in S$, 
\[-\frac{h+c\theta}{\alpha+\theta}\leq \Delta v^*(i)<0.\]    
\end{proof}
The ascent direction in (\ref{NL1}) is $(-\Delta v^*(i-1),-1)$. Therefore, we immediately have the following proposition.
\begin{Proposition}\label{PNL}
   For the infinite-horizon discounted reward optimality, $\mu$ corresponding to $a_i^*$, for all $i\geq 1$ should satisfy
    \[\partial _-f^*(\mu)\leq\frac{h+c\theta}{\alpha+\theta},\]
    where $\partial_-f^*$ denotes the left derivative of the function $f^*$.
        Furthermore, if  $f\geq 0$, then $\mu$   corresponding to $a_i^*$, for all $i\geq 1$ should satisfy
    \[\partial_+f^*(\mu)>0,\]
    where $\partial_+f^*$ denotes the right derivative of the function $f^*$.
\end{Proposition}
\begin{proof}
    The solution to the optimization problem in (\ref{NL1}), for all $i\geq1$, is given as
    \[\mu^*\in\{\mu:\partial_-f^*(\mu)\leq -\Delta v^*(i-1)\leq \partial_+f^*(\mu)\}.\]
We know from Proposition \ref{P31} that for all $i\geq1$, $\Delta v^*(i-1)\geq-\frac{h+c\theta}{\alpha+\theta}$ and the first result follows. If $f\geq 0$, then $-\frac{h+c\theta}{\alpha+\theta}\leq \Delta v^*(i)<0.$ Plugging this into the solution set gives the second inequality.
\end{proof}
The following two corollaries establish the optimal service rate for two special cases.

\begin{corollary}
    If the lower boundary $f^*$ is a linear function in $[\mu_-,\mu_+]$, then the optimal service rate in all states is either $\mu_-$ or $\mu_+$.
\end{corollary}
\begin{corollary}
    If $f\geq 0$  and the lower boundary  $f^*$ is non-increasing in $[\mu_-,\mu_+]$, then the optimal policy $\pi^*$ is a static policy which chooses service rate $\mu_+$ in all states.
\end{corollary}
Proposition \ref{PNL} further narrows the set of candidate optimal actions established in Theorem \ref{T34CH}. We next provide an explicit example to illustrate this. Suppose the projection of action space $A$ is a closed ball as shown in the gray area in Figure 3.1. Then, from Theorem \ref{T34CH}, we know the optimal actions are on the lower boundary $f^*$ which is the red and blue arc of the circle. Moreover, from Proposition \ref{PNL} we know that for the optimal actions, $\partial_-f^*(\mu)$ is bounded above by $\frac{h+c\theta}{\alpha+\theta}$ and below by $0$. Therefore, we only need to consider the blue arc of the circle when we search for  optimal actions.
\begin{figure}[H]
\centering
\begin{tikzpicture}[scale=1.2]

  \draw[->] (0,0) -- (4.2,0) node[right] {\footnotesize$\mu$};
  \draw[->] (0,0) -- (0,4.2) node[above] {\footnotesize$f$};
  \fill[gray!50, opacity=0.7, draw=black, thick] (2,2) circle (1);
  \draw[red, ultra thick] (1,2) arc[start angle=-180, end angle=0, radius=1];
  \draw[blue, ultra thick] (2,1) arc[start angle=-90, end angle=-10, radius=1];
  \draw[dashed, thick] (1.5,1) -- (2.5,1)
   node[below left] {\footnotesize slope $=0$};
  
  \draw[dashed, thick] (2.8980,1.3340) -- (3.0716,2.3188)
  node[right] {\footnotesize slope $=\frac{h+c\theta}{\alpha+\theta}$};

  \draw[->, thick] (2.8,0.9) -- (3.4,0.5)
  node[right] {\footnotesize possible ascent direction};
  
\end{tikzpicture}
\caption{Set of possible optimal actions when the projection of  $A$ is a closed ball}
\end{figure}

\subsection{Long-run average reward optimality}\label{S33}
In this section, we focus on characterizing the  policy that maximizes the long-run average reward. For any policy $\pi\in\Pi$, the long-run average reward when the initial state is $i$ is defined as 
\[g(i,\pi):=\liminf_{T\to\infty}\frac{1}{T}\mathbb{E}_i^\pi[\int_0^T R(t,X^\pi(t),\pi)\mathrm{d}t],\]
\noindent where the expectation is taken with respect to the continuous time Markov chain under policy 
$\pi$. The optimal gain is then defined as
\[g^*(i):=\sup_{\pi\in\Pi}g(i,\pi).\]

Because $b_2 = -\theta < 0$ in Theorem \ref{A1} and Theorem \ref{A2}(i) holds, $g^*(i) < +\infty$ for all $i \in S$.

We need the following result to introduce the optimality equations for the long-run average reward criterion.
\begin{theorem}\label{A4}
     For the long-run average reward model, if the Theorems \ref{A1}, \ref{A2} and \ref{A3} hold, and for each $\pi\in\mathbb{F}$,\\
     (i) the Markov process $\{X^\pi(t):t\geq0\}$ is irreducible with $P_\pi$ as the unique invariant probability measure under policy $\pi$,\\
(ii)   
    \[\sum_{j\geq k}q_{ij}(\pi)\leq \sum_{j\geq k}q_{i+1,j}(\pi) \mbox{ for all }i,k\in S,\ k\neq i+1,\]
(iii) for each $j>i_1>0$, there exist nonzero distinct states $i_1,i_2,\dots,i_m\geq j$ such that
    \[q_{i_1,i_2}(\pi)q_{i_2,i_3}(\pi)\dots q_{i_{m-1},i_m}(\pi)>0,\]
then, \\
(a) there exists a solution $(g^*,u)\in\mathbb{R}\times\mathbb{B}_\omega$ to the optimality equations
\begin{equation}\label{OE3}
    g^*(i)=\sup_{a\in A(i)}\{r(i,a)+\sum_{j\in S}q_{ij}(a)u(j)\}\qquad \forall\  i\in S
\end{equation}
where $u$ is unique up to additive constant vectors,\\
(b) the long-run average reward vector is a constant vector, i.e., $g^*(i)=g^*$, for all $i\in S$, and\\
(c) there exists a stationary deterministic policy $\pi \in\mathbb{F}$ that attains the maximal long-run average profit.
\end{theorem}
We know that in our model, the continuous time Markov chain is a birth-death process under any policy $\pi\in\mathbb{F}$. Therefore, conditions (i), (ii) and (iii) of Theorem \ref{A4}  are satisfied and there exists an optimal stationary deterministic policy $\pi^*\in\mathbb{F}$ that satisfies the optimality equations in (\ref{OE3}).
Next, we state the main result of this subsection. With a slight abuse of notation, we again use $\pi^*$ to denote the optimal policy and $a_i^*$ to denote the optimal action in state $i$ for the long-run average reward criterion.
\begin{theorem}\label{T36}
    The optimal service rate control policy $\pi^*$ that maximizes the long-run average profit is monotone non-decreasing in the number of customers in the system.
\end{theorem}

For the rest of this subsection, we first prove Theorem \ref{T36} and then explore the properties of the optimal service rate control policy. Note that the  optimality equations in ($\ref{OE3}$) are equal to 
\begin{equation}\label{EEv2}
    g^*=\left\{
\begin{array}{rcl}
    & \lambda r+\lambda u(1)-\lambda u(0) & \mbox{for }i=0,\\
    & \lambda r-hi-c(i-1)\theta-f(a_i^*)+\lambda u(i+1)\\
    &+(\mu(a_i^*)+(i-1)\theta)u(i-1)-(\lambda+\mu(a_i^*)+(i-1)\theta)u(i)&\mbox{for }i>0.\\
\end{array}
\right.
\end{equation}\par
\normalsize
From (\ref{EEv2}) with some algebra we have,
\normalsize
\begin{equation}\label{EEv3}
     0=\left\{
\begin{array}{rcl}
     &-h-f(a_1^*)+\lambda\Delta  u(1)-(\lambda+\mu(a_1^*))\Delta  u(0)& \mbox{for } i=0,  \\
     & -h-c\theta+\lambda\Delta  u(i+1)-(\lambda+i\theta)\Delta  u(i)+(i-1)\theta\Delta  u(i-1)\\
     &-f(a_{i+1}^*)+f(a_i^*)-\mu(a_{i+1}^*)\Delta  u(i)+\mu(a_i^*)\Delta u(i-1) & \mbox{for } i>0.\\
     
\end{array}
\right.
\end{equation}
Define $w^\prime(i)=-f(a_{i+1}^*)+f(a_i^*)-\mu(a_{i+1}^*)\Delta  u(i)+\mu(a_i^*)\Delta u(i-1)$. Then from (\ref{EEv3}), with some algebra we have,
\begin{equation}\label{OELR2}
    0=\left\{
\begin{array}{rcl}
    & f(a_1^*)-c\theta+\lambda\Delta^2 u(1)-(\lambda+\theta)\Delta^2 u(0)  \\
    &+(\mu(a_1^*)-\theta)\Delta u(0)+w^\prime(1)&\mbox{for } i=0,\\
    & \lambda\Delta^2 u(i+1)-(\lambda+(i+1)\theta)\Delta^2 u(i)+(i-1)\theta\Delta^2 u(i-1)\\
    & +w^\prime(i+1)-w^\prime(i)&\mbox{ for } i>0.
\end{array}
\right.
\end{equation}
\par

Similar to the infinite-horizon discounted reward criterion, we first outline the key steps and technical arguments required for the proof of Theorem \ref{T36}.
\begin{enumerate}
    \item We show that $\Delta u(i-1)$, together with the $\mu$ and $f$ functions, completely characterizes the optimal service rate control policy in state $i$ (Lemma \ref{L31}). Consequently, $\mu(a_i^*)$ is non-decreasing in $i$ if and only if $\Delta^2 u(i)$ is less than or equal to $0$. 
    \item We show that for $i\geq1$, $w^\prime (i)=-\gamma^\prime_{i-1}\Delta^2 u(i-1)$ for some constant vector $\{\gamma^\prime_i\}_{i=0}^\infty$ (Lemma \ref{L38C}). As a result, we can rewrite equation (\ref{OELR2}) as a system of equations for $\Delta^2 u(i),\ i\geq0$.
    \item We truncate this system of equations and consider the first $L$ number of equations (equation (\ref{EqLR})). We show that the solution to this truncated system of equations $\Delta^2 u_L(i), 0\leq i\leq L$ is less than or equal to $0$ (Lemmas \ref{L39C} and  \ref{L39}).
    \item We show that $\Delta^2 u_L(i)$ is monotone non-increasing in $L$ and pointwise converges to $\Delta^2 u(i)$. As a result, the solution to the original system of equations is less than or equal to $0$ (Lemma \ref{L39}).
    \item We conclude that the optimal service rate is monotone non-decreasing in the number of customers in the system.
\end{enumerate}

The following result which is similar to   Lemma \ref{Lemma31} characterizes the set of optimal actions in each state $i\geq1$.

\begin{lemma}\label{L31}
    Under long-run average reward optimality criterion, for $i\geq 1$, $a_i^*$   satisfies
    \[a_i^*\in\mathop{argmax}_{a\in A}\{-f(a)-\mu(a)\Delta u(i-1)\}.\]
\end{lemma}
Proof of the following lemma is similar to the proof of Lemma \ref{L32} and omitted.
\begin{lemma}\label{L38C}
    For $i\geq1 $, we have 
    \begin{equation}
        -\mu(a_i^*)\Delta^2u(i-1)\leq w^\prime(i)\leq -\mu(a_{i+1}^*)\Delta^2 u(i-1).
    \end{equation}
\end{lemma}
Therefore, there exists $\gamma^\prime_{i-1}\in[\min\{\mu(a_{i}^*),\mu(a_{i+1}^*)\},\max\{\mu(a_{i}^*),\mu(a_{i+1}^*)\}]$, such that $w^\prime(i)=-\gamma^\prime_{i-1}\Delta^2 u(i-1)$.
Let $G^\prime=\frac{f(a_1^*)-c\theta+(\mu(a_1^*)-\theta)\Delta u^*(0)}{\lambda}\leq0$ and construct a system of equations for a truncated system with $L\geq 1$ as
\small
\begin{equation}\label{EqLR}
    (\lambda+(i+1)\theta+\gamma_i^\prime)\Delta^2 u_L(i)=\left\{
\begin{array}{rcl}
     &\lambda G^\prime+\lambda\Delta^2 u_L(1)& \mbox{for } i=0,  \\
     & \lambda\Delta^2 u_L(i+1)+((i-1)\theta+\gamma_{i-1}^\prime)\Delta^2 u_L(i-1) & \mbox{for } 0<i\leq L,\\
     &0&\mbox{for } i>L.
     
\end{array}
\right.
\end{equation}
\normalsize
For $i\geq0$, define $p_i^\prime=\frac{(i+1)\theta+\gamma_i^\prime}{\lambda}$. Let $q_0^\prime=p_0^\prime+1$, $q_1^\prime=p_1^\prime q_0^\prime+q_0^\prime-p_0^\prime+\frac
{\theta}{\lambda}$ and $q_{i+2}^\prime=q_{i+1}^\prime p_{i+2}^\prime+q_{i+1}^\prime-p_{i+1}^\prime q_i^\prime+q_i^
\prime\frac{\theta}{\lambda}$ for $i\geq 0$. Note that $p_i$ and $q_i$ evaluated at $\alpha=0$ yield $p_i^\prime$ and $q_i^\prime$, respectively. We then immediately have the following lemma which is similar to Lemma \ref{L34}.
\begin{lemma}\label{L38p}
For $i\geq0$,
\[0<q_i^\prime p_{i+1}^\prime<q_{i+1}^\prime<q_i^\prime(p_{i+1}^\prime+1).\]
\end{lemma}
Then the following lemma follows as in Section \ref{S31}.
\begin{lemma}
     For $i\geq0$,
    \[\Delta^2 u(i+1)=q_i^\prime\Delta^2 u(0)-\Pi_{j=1}^{i}(p_j^\prime+1)G^
    \prime\]
    and for $0\leq i\leq L$
    \[\Delta^2 u_L(i+1)=q_i^\prime\Delta^2 u_L(0)-\Pi_{j=1}^{i}(p_j^\prime+1)G^\prime\]
    with the convention that $\Pi_{j=1}^0(p_j+1)=1$.
    
\end{lemma}
Following two lemmas are similar to Lemmas \ref{L37} and \ref{L36}, respectively.
 \begin{lemma}\label{L39C}
     Let $\Delta^2 u_L(i),i\in S$ be the unique solution to the system of equations in (\ref{EqLR}). If $G^\prime=0$, then $\Delta^2 u_L(i)=0$ for all $i\in S$ and $L\geq1$.
 \end{lemma}
\begin{lemma}\label{L39}
     Let $\Delta^2 u_L(i),i\in S$ be the unique solution to the system of equations in (\ref{EqLR}). If $G^\prime<0$, then\\
    (i) $\Delta^2 u_L(i)<0$ for all $i\leq L$;\\
    (ii) $\Delta^2u_L(i)>\Delta^2u_{L+1}(i)$ for all $i\leq L$.
    
\end{lemma}
 From Theorem \ref{A4}, we know that the bias vector $u$ is unique up to additive constant vectors. Therefore, $\Delta u$ and $\Delta^2 u$ are unique vectors. As a result, for any bias function $u$, $\Delta^2 u(i)<0$ for all $i\in S$ and the structure of the optimal policy remains the same under the long-run average reward optimality criterion. 
 This concludes the proof of Theorem \ref{T36}.
 Note that Theorem \ref{T36}
can also be shown by letting $\alpha\to0$ in the discounted model and proving that $v^*(i)$ of the discounted problem converges to real bias $u(i)$ of the long-run average reward problem (as shown in Theorem 4.1 of Guo et al. \cite{survey}) but we preferred to provide a more direct proof.\par
Furthermore, Theorem 4.1 of Guo et al. \cite{survey} also shows that for any policy $\pi\in\mathbb{F}$ and any sequence $\{\alpha_n: n\geq 1\}$ such that $\alpha_n \downarrow 0$ as $n\rightarrow \infty$, we have
\[\lim_{n\to\infty}\alpha_nv(i,\pi,\alpha_n)= g(i,\pi)\quad\mbox{and}\quad \lim_{n\to\infty} \Delta v^*(i,\alpha_n)= \Delta u(i).\]
Therefore, the properties of $g^*$ and the long-run average optimal policy $\pi^*$ will be similar to those of their counter parts in the discounted case.
\begin{Proposition}
    The  long-run average reward optimal policy does not depend on $r$.
\end{Proposition}
For $i\geq1$, optimal actions in state $i$ can again be determined by solving the following optimization problem
\begin{equation}
    \max_{(x,y)\in\{(\mu(a),f(a)),a\in A\}}\qquad -y-x\Delta u(i-1).\nonumber
\end{equation}
As a result, Theorem \ref{T34CH} holds for the long-run average reward criterion and the optimal actions can only be chosen on the lower boundary of the convex hull of the action space. The optimization problem is transformed into 
\begin{equation}\label{NL2}
    \max_{\mu\in[\mu_-,\mu_+]}\qquad -f^*(\mu)-\mu\Delta u(i-1).
\end{equation}
The next proposition provides bounds for the optimal gain $g^*$ and $\Delta u(i)$.
\begin{Proposition}\label{P35}
     For the long-run average reward criterion,  we have
\[g^*\geq \lambda(r-\frac{h+c\theta}{\theta})\quad\mbox{and}\quad \Delta u(i)\geq -\frac{h+c\theta}{\theta},\mbox{ } \forall i\in S.\]
Furthermore, if  $f\geq0$, we have
    \[\lambda(r-\frac{h+c\theta}{\theta})\leq g^*<\lambda r\mbox{ and} -\frac{h+c\theta}{\theta}\leq \Delta u(i)<0,\mbox{ } \forall i\in S.\]
\end{Proposition}
As a result, Proposition \ref{PNL} holds for the long-run average reward model with $\alpha=0$. Recall that $f^*$ is the lower boundary of $conv(H)$.
\begin{Proposition}\label{PNL2}
   For the long-run average reward optimality, $\mu$ corresponding to $a_i^*$, for all $i\ge 1$ should satisfy
    \[\partial _-f^*(\mu)\leq\frac{h+c\theta}{\theta}\]
      where $\partial_-f^*$ denotes the left derivative of the function $f^*$. Furthermore, if the service cost function $f$ is non-negative, then $\mu$ corresponding to   $a_i^*$, for all $i\geq 1$ should satisfy
    \[\partial_+f^*(\mu)>0\]
     where $\partial_+f^*$ denotes the right derivative of the function $f^*$.
\end{Proposition}
\begin{remark} Throughout the paper, we assume that customers only abandon the system while waiting in the queue. However, if customers can also abandon during service—after an exponential duration with rate $\theta_s$, incurring an abandonment cost of $c_s$—the framework remains valid. In this case, the action space $H = \{(\mu(a), f(a)) : a \in A\}$ is simply translated to the set $H + (\theta_s, c_s\theta_s)$, where the addition denotes the Minkowski sum of two sets (see, e.g., Theorem 3.1 in \cite{convex}). Under this extension, the  assumptions $f(\text{Idle}) = c\theta$ and $\mu(\text{Idle}) = \theta$ no longer apply; instead, they are replaced by $f(\text{Idle}) = \mu(\text{Idle}) = 0$.
\end{remark}
\section{Optimal policy when customers pay at service completion} \label{S4}

In this section, we consider the version of our problem where the reward $r$ is collected upon service completion rather than upon arrival. We begin by establishing that this formulation is mathematically equivalent to the original model. Moving forward, $a_i^{**}$ will denote the optimal action in state $i$ under both the infinite-horizon discounted and long-run average reward optimality criteria.
\begin{theorem}\label{T41}
    The version where customers pay the reward $r$ at the time of service completion  is equivalent to the original version where customers pay a reward of zero at the time of arrival and service cost function $F(a)=f(a)-\mu(a)r$.
\end{theorem}
\begin{proof}
   Note that when customers pay at the time of service completion, $q_{ij}(a)$ for all $a\in A$ remains the same.
    However, the reward function $r^\prime(i,a)$ is different. In particular, for $i\geq1$ and $a\in \mathcal{A}$, 
   \[r^\prime(i,a)=\mu(a)r-hi-c(i-1)\theta-f(a)=-hi-c(i-1)\theta-[f(a)-\mu(a)r]\]
    and for $i\geq 0$ and $a=Idle$, we have
\[r^\prime(i,Idle)=-hi-c\theta i.\]
Hence, the result follows.
\end{proof}
For the long-run average reward optimality, we provide another equivalence between the two versions in the Appendix.

When the customers pay the reward $r$ at the time of service completion, for a policy $\pi\in\Pi$, define the reward function as
\[R^\prime(t,i,\pi)=\int_{A(i)}r^\prime(i,a)\rho_t(\mathrm{d}a|i).\]\par
Recall that $\rho_t(\mathrm{d}a|i)$ is the probability that the policy chooses actions in a Borel set $B\in\mathcal{B}(A(i))$ at time $t$. For the infinite-horizon discounted reward criterion with discount factor $\alpha>0$, if the system starts in state $i\in S$, then the value function of policy $\pi$ is defined as 
\[V(i,\pi,\alpha)=\mathbb{E}_i^\pi[\int_0^\infty e^{-\alpha t}R^\prime(t,i,\pi)\mathrm{d}t]\]
and the optimal discounted reward is equal to
\[V^*(i,\alpha)=\sup_{\pi\in\Pi}V(i,\pi,\alpha).\]\par
In order to simplify the notation, we will write $V^*(i)$ instead of $V^*(i,\alpha)$ for the optimal value function. Then we have the following optimality equation
\begin{equation}\label{OEModel2}
   \alpha V^*(i)=\sup_{a\in A(i)}\{r^\prime(i,a)+\sum_{j\in S}q_{ij}(a)V^*(j)\}\qquad \forall\  i\in S.
\end{equation}
Similarly, long-run average reward under policy $\pi$ when the initial state is $i\in S$ is defined as
\[\mathcal{G}(i,\pi)=\liminf_{T\to\infty}\frac{1}{T}\mathbb{E}_i^\pi[\int_0^T R^\prime(t,i,\pi)\mathrm{d}t]\]
and the optimal gain $\mathcal{G}^*$ and bias vector $U$ satisfy
\[\mathcal{G}^*(i):=\sup_{\pi\in\Pi}\mathcal{G}(i,\pi) \] and
\[\mathcal{G}^*(i)=r^\prime(i,a_i^{**})+\sum_{j\in S}q_{ij}(a_i^{**})U(j) \qquad\forall i\in S,\]
where the second equality follows from Theorem \ref{A4}. Note that $\mathcal{G}^*(i)=\mathcal{G}^*$ for all $i\in S$  and $U$ is unique up to constant vectors (see Theorem \ref{A4}).\par
The new service  cost function $F(a)=f(a)-\mu(a)r$ is still continuous in $a\in A$. Furthermore, Theorem \ref{T41} does not depend on the optimality criterion. Therefore,  all the theorems and lemmas of Section \ref{S3} apply to the version where customers pay at the time of service completion. Then from Theorems \ref{T33} and \ref{T36}, we have the following results.
\begin{corollary}
    For $i\in S$, both $\Delta V^*(i)$ and $\Delta U(i)$ are non-increasing in $i$. For the version where customers pay at the time of service completion, the optimal service rate control policy that maximizes the infinite-horizon discounted reward (long-run average reward) is monotone non-decreasing in state $i\in S$.
\end{corollary}
From Proposition \ref{P31}, we know that the optimal policy does not depend on the reward paid at the time of customer arrival. However, the reward paid after service completion will affect the optimal policy. Next, we explore how the optimal service rate  changes with respect to the reward.\par
First, we consider the infinite-horizon discounted reward criterion. Expanding the optimality equation in (\ref{OEModel2}), we have
\small
\begin{equation}\label{OE3model2}
     \alpha\Delta V^*(i)=\left\{
\begin{array}{rcl}
     &-h-f(a_1^{**})+\mu(a_1^{**})r+\lambda\Delta  V^*(1)-(\lambda+\mu(a_1^{**}))\Delta  V^*(0)& \mbox{for } i=0,  \\
     & -h-c\theta+\lambda\Delta  V^*(i+1)-(\lambda+i\theta)\Delta  V^*(i)+(i-1)\theta\Delta  V^*(i-1)\\
     &-f(a_{i+1}^{**})+f(a_i^{**})+(\mu(a_{i+1}^{**})-\mu(a_i^{**}))r\\
     &-\mu(a_{i+1}^{**})\Delta  V^*(i)+\mu(a_i^{**})\Delta V^*(i-1) & \mbox{for } i>0.\\
     
\end{array}
\right.
\end{equation}
\normalsize

We next consider two systems with exactly the same system parameters except for the reward paid after service completion. Assume that customers pay $r_1$ and $r_2$ at the time of service completions in systems one and two, respectively. Let $V_1^*(i)$ and $V_2^*(i)$ be the corresponding optimal value functions and $a_{i,1}^{**}$ and $a_{i,2}^{**}$ be the corresponding optimal actions in state $i\in S$ for the two systems. Without loss of generality, suppose $r_1<r_2$ and define $\Delta r=r_2-r_1>0$ and $\Delta^2 \mathcal{V}_r(i)=\Delta V_2^*(i)-\Delta r-\Delta V_1^*(i)$. We first prove the following lemma.
\begin{lemma}\label{L41}
    For any $i\in S$,
    \begin{eqnarray}
    -\mu(a_{i,1}^{**})\Delta^2 \mathcal{V}_r(i-1)\leq&\mu(a_{i,1}^{**})\Delta V_1^*(i-1)+f(a_{i,1}^{**})-\mu(a_{i,1}^{**})r_1&\nonumber\\
    &-\mu(a_{i,2}^{**})\Delta V_2^*(i-1)-f(a_{i,2}^{**})+\mu(a_{i,2}^{**}) r_2&\leq-\mu(a_{i,2}^{**})\Delta^2 \mathcal{V}_r(i-1).\nonumber
    \end{eqnarray}
\end{lemma}
\begin{proof}
    From the optimality equation in (\ref{OEModel2}) we know that
    \[a_{i,1}^{**}\in\mathop{argmax}_{a\in A}\{-f(a)+\mu(a)r_1-\mu(a)\Delta V_1(i-1)\}\]
    and
    \[a_{i,2}^{**}\in\mathop{argmax}_{a\in A}\{-f(a)+\mu(a)r_2-\mu(a)\Delta V_2(i-1)\}.\]
    Therefore, 
    \[\mu(a_{i,1}^{**})\Delta V_1^*(i-1)+f(a_{i,1}^*)-\mu(a_{i,1}^{**})r_1\leq\mu(a_{i,2}^{**})\Delta V_1^*(i-1)+f(a_{i,2}^{**})-\mu(a_{i,2}^{**})r_1\]
    and
    \[-\mu(a_{i,2}^{**})\Delta V_2^*(i-1)-f(a_{i,2}^{**})+\mu(a_{i,2}^{**})r_2\geq-\mu(a_{i,1}^{**})\Delta V_2^*(i-1)-f(a_{i,1}^{**})+\mu(a_{i,1}^{**})r_2.\]
    Hence, the result follows.
\end{proof}
Lemma \ref{L41} is similar to Lemma \ref{L32} and we immediately have the following corollaries.
\begin{corollary}\label{C42}
    If for $i\geq 1$, $\Delta^2 \mathcal{V}_r(i-1)<0$, then $\mu(a_{i,1}^{**})\leq \mu(a_{i,2}^{**})$.
\end{corollary}
The following corollary again follows from the Intermediate Value Theorem (see Theorem 4.23 in Rudin \cite{rudin}).
\begin{corollary}\label{C43}
    There exist constants $\tau_i\in[\min\{\mu(a_{i,1}^{**}),\mu(a_{i,2}^{**})\},\max\{\mu(a_{i,1}^{**}),\mu(a_{i,2}^{**})\}]$, for $i\geq1$ such that 
    \begin{equation}\label{E43}
    \mu(a_{i,1}^{**})\Delta V_1^*(i-1)+f(a_{i,1}^{**})-\mu(a_{i,1}^{**})r_1
    -\mu(a_{i,2}^{**})\Delta V_2^*(i-1)-f(a_{i,2}^{**})+\mu(a_{i,2}^{**}) r_2=-\tau_i\Delta^2\mathcal{V}_r(i-1).
    \end{equation}
\end{corollary}
Subtracting (\ref{OE3model2}) for $V^*_1$ from (\ref{OE3model2}) for $V^*_2$ and plugging  (\ref{E43}) into the system of equations, we have

\begin{equation}\label{OE3model25}
     \alpha\Delta^2\mathcal{V}_r(i)+\alpha \Delta r=\left\{
\begin{array}{rcl}
     &\lambda\Delta^2\mathcal{V}_r(1)-(\lambda+\tau_1)\Delta^2\mathcal{V}_r(0)& \mbox{for } i=0,  \\
     & \lambda\Delta^2\mathcal{V}_r(i+1)-(\lambda+i\theta+\tau_{i+1})\Delta^2\mathcal{V}_r(i)&\\
     &+((i-1)\theta+\tau_i)\Delta^2\mathcal{V}_r(i-1)-\theta \Delta r & \mbox{for } i>0.\\
     
\end{array}
\right.
\end{equation}
\normalsize
\begin{theorem}\label{T43}
For all $i\in S$,
    \[
    \Delta^2\mathcal{V}_r(i)<0.
    \]
\end{theorem}
\begin{proof}
    Rewriting the equations in (\ref{OE3model25}), we have
\begin{equation}\label{E431}(\alpha+\lambda+\tau_1)\Delta^2\mathcal{V}_r(0)+\alpha \Delta r=\lambda\Delta^2\mathcal{V}_r(1),
\end{equation}
\begin{equation}\label{E432}
(\alpha+\lambda+i\theta+\tau_{i+1})\Delta^2\mathcal{V}_r(i)+(\alpha+\theta)\Delta r=\lambda\Delta^2\mathcal{V}_r(i+1)+((i-1)\theta+\tau_{i})\Delta^2\mathcal{V}_r(i-1) \qquad \mbox{for } i>0.   
\end{equation}
\normalsize
Let $k\geq2$. Then summing up the equations (\ref{E431}) and (\ref{E432}) for $i=1,\dots,k-1$, we have
\begin{equation}\label{E433}
\lambda\Delta^2\mathcal{V}_r(k)=\lambda\Delta^2 \mathcal{V}_r(0)+\sum_{i=1}^k\alpha\Delta^2 \mathcal{V}_r(i-1)+((k-1)\theta+\tau_k)\Delta^2 \mathcal{V}_r(k-1)+k\alpha \Delta r+(k-1)\theta \Delta r.
\end{equation}
It follows from (\ref{E431}) that, if $\Delta^2 \mathcal{V}_r(0)\geq0$,  $\lambda\Delta^2 \mathcal{V}_r(1)\geq \alpha \Delta r$. Then, from (\ref{E433}), we have $\lambda\Delta^2\mathcal{V}_r(i)\geq (i-1)\theta \Delta r$ for all $i\geq 1$ and the sequence $\{\Delta^2 \mathcal{V}_r(i)\}_{i\in S}$ is not bounded above which is a contradiction since $V^*(i)\in\mathbb{B}_\omega$ (see Section \ref{S31}). Hence, we must have $\Delta^2 \mathcal{V}_r(0)<0$. \par
We now use induction. Suppose $\Delta^2\mathcal{V}_r(i)<0$ for all $i<k$, we will show that $\Delta^2\mathcal{V}_r(k)<0$. If $\Delta^2\mathcal{V}_r(k)\geq0$, then from (\ref{E433}) and the induction hypothesis, we have
\[0\leq\lambda\Delta^2\mathcal{V}_r(k)<\lambda\Delta^2\mathcal{V}_r(0)+\sum_{i=1}^k\alpha\Delta^2\mathcal{V}_r(i-1)+k\alpha \Delta r+(k-1)\theta \Delta  r. \]
Therefore,

\begin{eqnarray}
\lambda\Delta^2\mathcal{V}_r(k+1)&=&\lambda\Delta^2 \mathcal{V}_r(0)+\sum_{i=1}^k\alpha\Delta^2\mathcal{V}_r(i-1)+\alpha\Delta^2\mathcal{V}_r(k)+(k\theta+\tau_{k+1})\Delta^2 \mathcal{V}_r(k)+(k+1)\alpha \Delta r +k\theta \Delta r\nonumber\\
    &>&\alpha\Delta^2\mathcal{V}_r(k)+(k\theta+\tau_{k+1})\Delta^2 \mathcal{V}_r(k)+\alpha \Delta r+\theta \Delta r\nonumber\\
    &>&\theta \Delta r.\nonumber
\end{eqnarray}
Similarly, for all $j>k+1$, we have
\begin{eqnarray}
    \lambda\Delta^2\mathcal{V}_r(j)&>&\sum_{i=k+1}^j \alpha\Delta^2\mathcal{V}_r(i-1)+((j-1)\theta+\tau_j)\Delta^2\mathcal{V}_r(j-1)+(j-k)\alpha \Delta r+(j-k)\theta \Delta r\nonumber\\
    &>&(j-k)\theta \Delta r,\nonumber
\end{eqnarray}
which again implies that $\Delta^2\mathcal{V}_r(i),i\in S$ is not bounded above, which is a contradiction. Consequently, we must have $\Delta^2 \mathcal{V}_r(k)<0$. This completes the proof that $\Delta^2\mathcal{V}_r(i)<0$ for all $i\in S$.
\end{proof}
When customers pay at the time of service completion, we have from (\ref{OEModel2}) that for $i\geq 1$
\begin{equation}\label{NLModel2}
    \max_{(x,y)\in\{(\mu(a),f(a)),a\in A\}}\qquad -y-x\left(\Delta V^*(i-1)-r\right).
\end{equation}
As a result, we know that Theorem \ref{T34CH} also holds for this version. In order to obtain similar bounds (as in Proposition \ref{PNL}) for the optimal actions, we first need to obtain bounds for $V^*$ and $\Delta V^*$.
\begin{Proposition}\label{P41}
For the infinite-horizon discounted reward criterion, for all $i\in S$, we have
\[V^*(i)\geq \frac{\lambda}{\alpha}(r-\frac{h+c\theta}{\alpha+\theta})-\frac{h+c\theta}{\alpha+\theta}i\quad \mbox{and}\quad\Delta V^*(i)\geq -\frac{h+c\theta}{\alpha+\theta}.\]
Furthermore, if  $f\geq0$, we have
\[-\frac{\lambda}{\alpha}\frac{h+c\theta}{\alpha+\theta}\leq V^*(0)<\frac{\lambda}{\alpha}r\mbox{ and }-\frac{h+c\theta}{\alpha+\theta}\leq \Delta V^*(i)<r \mbox{ for all $i\in S$}.\]    
\end{Proposition}
Proof of Proposition \ref{P41} is similar to the proof of  Proposition \ref{P31}. 
Note that the ascent direction of the optimization problem in (\ref{NLModel2}) is $(-\Delta V^*(i-1)+r,-1)$. From Proposition \ref{P41}, we know that when the cost function $f$ is non-negative, 
\[0<-\Delta V^*(i-1)+r\leq \frac{h+c\theta}{\alpha+\theta}+r.\]
Therefore, we immediately have the following proposition.
\begin{Proposition}\label{PNL3}
     For the infinite-horizon discounted reward criterion, let $f^*$ be the lower boundary of $conv(H)$. Then, $\mu$ corresponding to $a_i^{**}$, for all $i\geq 1$, should satisfy
    \[\partial _-f^*(\mu)\leq\frac{h+c\theta}{\alpha+\theta}+r\]
where $\partial_-f^*$ denotes the left derivative of the function $f^*$.    Furthermore, if $f\geq 0$, then $\mu$ corresponding to $a_i^{**}$, for all $i\geq 1$, should satisfy
    \[\partial_+f^*(\mu)>0\]
      where $\partial_+f^*$ denotes the right derivative of the function $f^*$.
\end{Proposition}
Note that Proposition \ref{PNL3} also illustrates that in the version where customers pay at the time of service completion, the set of possible optimal actions is still characterized  using the original service cost function $f$ rather than $F$.

From Theorem \ref{T43} we know that  if $r_1<r_2$, then $\Delta V_2^*(i)-r_2<\Delta V_1^*(i)-r_1$. As a result, in the optimization problem (\ref{NLModel2}), the coefficient of $\mu(a)$ increases as the reward $r$ increases, which leads to our main result of this section.

\begin{theorem}\label{T43r}
    When customers pay the reward $r$ at the time of service completion, under the infinite horizon discounted optimality criterion, optimal service rate $\mu(a_i^{**})$ for all $i\in S$ is non-decreasing in $r$. Furthermore, if $f\geq 0$, then $f(a_i^{**})$ for all $i\in S$ is also non-decreasing in $r$.
\end{theorem}
The following result compares the optimal actions of the two versions (payment at the time of arrival and service completion) and follows immediately from Proposition \ref{Pro31} and Theorem \ref{T43r}.
\begin{corollary}\label{C44}
   For the infinite-horizon discounted reward criterion, we have 
    \[\mu(a_i^*)\leq \mu(a_i^{**})\quad \forall \ i\geq1.\]
    Furthermore, if $f\geq 0$ , we have
    \[f(a_i^*)\leq f(a_i^{**})\quad \forall \ i\geq1.\]
\end{corollary}
Intuitively, if the system receives a reward when the customers complete their service, then the system is more willing to retain customers  and avoid abandonments as the reward increases. Consequently, the optimal service rate in any state tends to increase even though the service cost $f$ might also increase.\par

For the long-run average reward criterion, for $i
\in S$, define  $\Delta^2U_r(i)=\Delta U(i)-r-\Delta u(i)$. Since $\Delta u(i)$ and $\Delta U(i)$ are unique, $\Delta^2U_r(i)$ is well-defined. Using the arguments in the proof of Theorem \ref{T43} with $\alpha=0$, we can show that the results of Lemma \ref{L41} and Theorem \ref{T43} hold with 
$\Delta^2\mathcal{V}_r(i)$ replaced by $\Delta^2U_r(i)$. Thus, we have $\Delta^2U_r(i)<0$ for all $i\in S$. Note that as in Section \ref{S33}, we have for any sequence $\alpha_n$ such that $\alpha_n \downarrow 0$ as $n\rightarrow \infty$,
\[\lim_{n\to\infty}\alpha_nV(i,\pi,\alpha_n)= \mathcal{G}(i,\pi)\quad\mbox{and}\quad \lim_{n\to\infty} \Delta V^*(i,\alpha_n)= \Delta U(i).\]
Propositions \ref{newp1} and \ref{newp2} follow immediately  from Propositions \ref{P41} and \ref{PNL3}, respectively.

\begin{Proposition} \label{newp1}
  For the long-run average reward criterion,we have 
\[\mathcal{G}^*\geq \lambda(r-\frac{h+c\theta}{\theta})\quad\mbox{and for all $i\in S$,}\quad \Delta U^*(i)\geq -\frac{h+c\theta}{\theta}.\]
Furthermore, if $f\geq0$, we have
\[-\lambda\frac{h+c\theta}{\theta}\leq \mathcal{G}^*<\lambda r\quad\mbox{and for all $i\in S$,}
-\frac{h+c\theta}{\theta}\leq \Delta U^*(i)<r.\]    
\end{Proposition}
\begin{Proposition} \label{newp2}
     For the long-run average reward criterion, let $f^*$ be the lower boundary of $conv(H)$. Then, $\mu$ corresponding to $a_i^{**}$, for all $i\geq 1$, should satisfy
    \[\partial _-f^*(\mu)\leq\frac{h+c\theta}{\theta}+r\]
   where $\partial_-f^*$ denotes the left derivative of the function $f^*$.    Furthermore, if $f\geq 0$ is non-negative, then $\mu$ corresponding to $a_i^{**}$, for all $i\geq 1$, should satisfy
    \[\partial_+f^*(\mu)>0\]
     where $\partial_+f^*$ denotes the right derivative of the function $f^*$.
\end{Proposition}
Finally, Theorem \ref{T43r} and a more general version of Corollary \ref{C44} also hold for the long-run average reward criterion. 
\begin{theorem}\label{T45}
  When customers pay the reward $r$ at the time of service completion, under the long-run average optimality criterion, optimal service rate $\mu(a_i^{**})$ for all $i\in S$ is non-decreasing in $r$. Furthermore, if $f\geq 0$, then $f(a_i^{**})$ for all $i\in S$ is also non-decreasing in $r$.
\end{theorem}

\begin{corollary}\label{C45}
 We have 
    \[\mu(a_i^*)\leq \mu(a_i^{**})\quad \forall \ i\geq1\]
    and
    \[\mathcal{G}^*\leq g^*.\]
    Furthermore, if $f\geq 0$, we have
    \[f(a_i^*)\leq f(a_i^{**})\quad \forall \ i\geq1.\]
    \end{corollary}
Note that $\mathcal{G}^*\leq g^*$ since $\mathcal{G}(i,\pi)$ and $g(i,\pi)$ do not depend on the initial state $i$ and $\mathcal{G}(i,\pi)\leq g(i,\pi)$ for all $\pi\in\mathbb{F}$.

    Corollary \ref{C45} highlights a key distinction between systems with and without customer abandonments. Under the long-run average reward criterion without abandonments (i.e., $\theta = 0$), the optimal policy—provided it exists—is invariant to whether rewards are collected at arrival or upon service completion. Indeed, substituting $\alpha = \theta = 0$ into equation \eqref{OE3model25} yields $\Delta^2\mathcal{V}_r(i) = 0$ for all $i \in S$, demonstrating that the optimal policy remains unchanged under either reward structure. However, in the presence of abandonments, as stated in Corollary \ref{C45}, the timing of customer payments matters, since the optimal service rate when customers pay upon arrival is less than or equal to the optimal service rate when customers pay upon service completion.

\begin{remark}
   We note that a hybrid formulation where rewards are accrued both at arrival and upon service completion is also covered by this framework. In this case, for both optimality criteria, the optimal policy only depends on the reward paid at the time of service completion and hence, all results of this section hold.
\end{remark}

Corollaries \ref{C44} and \ref{C45} imply that the optimal service rate in each state is non-decreasing in the reward $r$. Applying the same proof technique, we immediately obtain the following proposition, which asserts the same monotonicity property with respect to the holding cost $h$ and abandonment cost $c$.
\begin{Proposition}
    For a state $i\in S$, the optimal service rate $\mu(a_i^*)$ and $\mu(a_i^{**})$ are both non-decreasing in $h$ and $c$. Furthermore, if $f\geq0$, then $f(a_i^*)$ and $f(a_i^{**})$ are also non-decreasing in $h$ and $c$.
\end{Proposition}
To summarize, we have shown that the optimal service rate in each state is non-decreasing in $r$, $c$, and $h$ under both optimality criteria. In contrast to these parameters, the optimal service rate is not necessarily non-decreasing in the abandonment rate $\theta$; we defer a discussion of this exception to Section \ref{S6}.

\section{Systems with finite buffer}\label{S5}
In this section, we focus on the version of the problem when the buffer size is finite. Theorem \ref{T51} states the interesting result that when there is limited buffer space,  optimal service rate policy is not necessarily monotone.

\begin{theorem}\label{T51}
    For the finite buffer model, the optimal service rate first increases and then decreases with the number of customers in the system under both optimality criteria.
\end{theorem}
We first prove Theorem \ref{T51} under the infinite-horizon discounted reward criterion, and then extend this result to the long-run average reward criterion.

We assume without loss of generality that $N > 2$. This is because the problem is trivial when $N = 1$, as there is only a single state where the service rate must be chosen, and when $N = 2$, the optimal service rate policy is inherently monotonic. Since $N < \infty$, an arriving customer finding $N$ in the system leaves. Define $v_{F}^*(i)$ as the optimal infinite horizon discounted reward when the initial state is $i$ and discount factor is $\alpha$. Let $a_{i,F}^*$ be the optimal action in state $i$, $i=0,\dots,N$ under this optimality criterion. Note that ``F" stands for the finite buffer system.  Then, the optimality equations for the infinite horizon discounted optimality criterion are
\[\alpha v_{F}^*(i)=\sup_{a\in A(i)}\{r(i,a)+\sum_{j=0}^{N}q_{ij}(a)v_{F}^*(j)\}\qquad  i=0,1,\dots,N-1,\]
and
\[\alpha v_{F}^*(N)=\sup_{a\in A}\{r(N,a)-\lambda r+\sum_{j=0}^Nq_{Nj}(a)v_{F}^*(j)\}.\]
As a result, $v_{F}^*(i)$ satisfies the optimality equations in (\ref{EEv1}) for $i=0,\dots,N-1$. When $i=N$, we have
\[\alpha v_{F}^*(N)=-Nh-c(N-1)\theta-f(a_{N,F}^*)+(\mu(a_{N,F}^*)+(N-1)\theta)v_{F}^*(N-1)-(\mu(a_{N,F}^*)+(N-1)\theta)v_{F}^*(N).\]\par
Then, $\Delta v_{F}^*(i)$ satisfies the equations in (\ref{OE311}) for $i=0,\dots,N-2$ and when $i=N-1$, we have
\[\alpha\Delta v_{F}^*(N-1)=-h-c\theta-\lambda r-(\lambda+(N-1)\theta)\Delta v_{F}^*(N-1)+(N-2)\theta\Delta v_{F}^*(N-2)+w_F(N-1)\]
where $w_F(i)=-f(a_{i+1,F}^*)+f(a_{i,F}^*)-\mu(a_{i+1,F}^*)\Delta  v_F^*(i)+\mu(a_{i,F}^*)\Delta v_F^*(i-1)$. \par
Similarly, we have $\Delta^2 v_{F}^*(i)$ satisfying the equations in (\ref{OE41}) for $i=0,\dots,N-3$ and when $i=N-2$,
\[(\alpha+\lambda+(N-1)\theta+\gamma_{N-2,F})\Delta^2 v_{F}^*(N-2)=-\lambda r-\lambda \Delta v_{F}^*(N-1)+((N-3)\theta+\gamma_{N-3,F})\Delta^2 v_{F}^*(N-3)\]
where $\gamma_{i,F}\in[\min\{\mu(a_{i,F}^*),\mu(a_{i+1,F}^*)\},\max\{\mu(a_{i,F}^*),\mu(a_{i+1,F}^*)\}]$ is again obtained using the Intermediate Value Theorem with $w_F(i)=-\gamma_{i-1,F}\Delta v_{F}^*(i-1)$ for $i=1,\dots,N-1$. Let $D=-r-\Delta v_{F}^*(N-1)$, then we have
\footnotesize
\begin{equation}\label{Finite} 
 (\alpha+\lambda+(i+1)\theta+\gamma_{i,F})\Delta^2 v_{F}^*(i)=\left\{
\begin{array}{rcl}
     &\lambda G+\lambda\Delta^2 v_{F}^*(1)& \mbox{for } i=0,  \\
     & \lambda\Delta^2 v_{F}^*(i+1)+((i-1)\theta+\gamma_{i-1,F})\Delta^2 v_{F}^*(i-1) & \mbox{for } 0<i< N-2,\\
     &\lambda D+((N-3)\theta+\gamma_{N-3,F})\Delta^2 v_{F}^*(N-3)&\mbox{for } i=N-2.
\end{array}
\right.
\end{equation}
\normalsize
Note that this system of equations is the same as the one in (\ref{EqL})  with $L=N-2$ and $\Delta^2 v_L(N-1)=D$. Therefore, with some algebra, we can compute $\Delta^2 v_{F}^*(i)$ as 
\begin{equation}\label{EF1}
    \Delta^2 v_{F}^*(0)=\frac{D}{q_{N-2,F}}+\Pi_{j=1}^{N-2}(p_j+1)\frac{G_F}{ q_{N-2,F}}
\end{equation} 
and 
\begin{equation}\label{EF2}
    \Delta^2 v_{F}^*(i)=\frac{q_{i-1,F}}{q_{N-2,F}}D+\frac{q_{i-1,F}}{q_{N-2,F}}\Pi_{j=1}^{N-2}(p_{j,F}+1)G_F-\Pi_{j=1}^{i-1}(p_{j,F}+1)G_F\qquad\mbox{for }i=1,\dots,N-2,
\end{equation} 
where $G_F$ is defined as in Section \ref{S31} with $\Delta v^*(0)$ replaced by $\Delta v_F^*(0)$ and $p_{i,F},q_{i,F}$ for $i=0,\dots,N-2$ are as defined in Section \ref{S31} with  $\gamma_i$ replaced by $\gamma_{i,F}$.\par
Note that since we do not know the sign of $D$, we can neither guarantee $\Delta^2v_F^*(i)<0$ nor the monotonicity of the optimal policy. However, the following lemma provides a relationship between $\Delta^2v_F^*(i)$ and $\Delta^2v_F^*(i+1)$ for $i=0,\dots,N-3$, which leads to a different structure of the optimal policy than the infinite buffer case.
\begin{lemma}\label{L51}
    For $i=0,\dots,N-3$, if $\Delta^2 v_{F}^*(i)>0$, then $\Delta^2 v_{F}^*(i+1)>0$.
\end{lemma}
\begin{proof}
    For $i=0$, if $\Delta^2 v_{F}^*(0)>0$, then as stated in Section \ref{S31}, it follows from Corollary \ref{C31} that $G_F\leq0$ and we have
    \[\Delta^2 v_{F}^*(1)=\frac{q_{0,F}}{q_{N-2}}D+\frac{q_{0,F}}{q_{N-2,F}}\Pi_{j=1}^{N-2}(p_{j,F}+1)G_F-G_F=q_{0,F}\Delta^2v_{F}^*(0)-G_F>q_{0,F}\Delta^2v_{F}^*(0).\]
      As in  Lemma \ref{L34}, for $i=0,\dots,N-3$, we have $q_{i+1,F}<q_{i,F}(p_{i+1,F}+1)$ and $q_{i,F}>0$. As a result, $\frac{q_{i,F}}{q_{i-1,F}}<p_{i,F}+1$ and
    \begin{eqnarray}
        \frac{q_{i,F}}{q_{i-1,F}}\Delta^2v_{F}^*(i)&=&\frac{q_{i,F}}{q_{N-2,F}}D+\frac{q_{i,F}}{q_{N-2,F}}\Pi_{j=1}^{N-2}(p_{j,F}+1)G_F-\frac{q_{i,F}}{q_{i-1,F}}\Pi_{j=1}^{i-1}(p_{j,F}+1)G_F\nonumber\\
        &\leq&\frac{q_{i,F}}{q_{N-2,F}}D+\frac{q_{i,F}}{q_{N-2,F}}\Pi_{j=1}^{N-2}(p_{j,F}+1)G_F-\Pi_{j=1}^{i}(p_{j,F}+1)G_F\nonumber\\
        &=&\Delta^2v_{F}^*(i+1).\nonumber
    \end{eqnarray}
    Therefore, if $\Delta^2 v_{F}^*(i)>0$, then $\Delta^2 v_{F}^*(i+1)>0$.
\end{proof}
\begin{lemma}
    For the finite buffer model with $N\geq 2$, there exists a state $1\leq k\leq N$ such that $\mu(a_i^*)$ is  non-decreasing in states $1\leq i\leq k$ and non-increasing in states $k\leq i\leq N$.
\end{lemma}
\begin{proof}
    If $\Delta^2 v_{F}^*(i)\leq0$ for all $i=0,\dots,N-2$, then we know that $\mu(a_i^*)$ is monotone non-decreasing in $i$ and $k=N$. Otherwise, define $k=\min\{0\leq i\leq N-2,\Delta^2 v_{F}^*(i)>0\}+1$, then from Lemma \ref{L51} we conclude that $\Delta^2 v_{F}^*(i)\leq0$ for $0\leq i< k-1$ and $\Delta^2 v_{F}^*(i)>0$ for  $k-1\leq i\leq N-2$. Consequently, $\mu(a_i^*)$ is  non-decreasing in states $1\leq i\leq k$ and  non-increasing in states $k\leq i\leq N$. 
\end{proof}
This concludes the proof of Theorem \ref{T51} under the infinite horizon discounted reward criterion. For the long-run average reward criterion, let $g_F^*$ be the optimal long-run average reward for the finite buffer model. Note again that $g_F^*$ does not depend on the initial state. Then, $g_F^*$ and the bias function $u_F(i)$ satisfy
\[g_F^*=\sup_{a\in A(i)}\{r(i,a)+\sum_{j=0}^Nq_{ij}(a)u_F(j)\}\qquad\forall i=0,\dots,N-1,\]
and
\[g_F^*=\sup_{a\in A}\{r(N,a)-\lambda r+\sum_{j=0}^Nq_{Nj}(a)u_F(j)\}.\]
In the system of equations in (\ref{Finite}), let $\alpha=0$ and substitute the bias function $u_F(i)$ for $v_F^*(i)$. Then we can show that Lemma \ref{L51} holds for $\Delta^2 u_F(i)$. As a result, Theorem \ref{T51} also holds under the long-run average reward criterion.\par 

Non-monotonicity of optimal service rates in a finite buffer queueing system without abandonments was also observed by Crabill \cite{crab1} who provided a numerical example illustrating this phenomenon. Intuitively, as the queue length approaches its upper bound, the influx of new customers slows down, meaning that the holding and abandonment costs per unit time no longer increase significantly. As a result, the optimal policy shifts toward a lower service rate to conserve costs, rather than aggressively clearing the system. For instance, consider a system with $\lambda=2,\ h=c=1.9,\ \theta=0.3$, and a buffer size of $10$. Suppose there are two possible service modes, $a_1$ and $a_2$, with $\mu(a_1)=1,\ f(a_1)=1$ and $\mu(a_2)=2,\ f(a_2)=8$. The optimal service rate control policy uses service mode $a_1$ in states $s=2,3,7,8,9,10$, and service mode $a_2$ in states $s=4,5,6$. Specifically, for the CTMDP in this model, $\Delta^2 v^*_F(i)$ is negative for the first five states and positive for states larger than $5$ (i.e., $k=6$ in Theorem \ref{T51}). Hence, it is optimal to use the lower service rate not only in the boundary state $N$ but also in the smaller states $7$, $8$, and $9$. Although the value of $k$ depends on the system parameters, this phenomenon observed in finite buffer systems appears robust across various parameter settings. In the next section, we show that this phenomenon vanishes as the buffer size tends to infinity.

\begin{remark}
 Note that the structure of the optimal policy in Theorem \ref{T51} remains the same if there is a rejection cost for the customers arriving when the system is full. Having a rejection cost only changes the value of $D$ which our proof does not depend on. 
\end{remark}

Finally, the equivalency in Theorem \ref{T41} also holds for the finite buffer case. Consequently, Theorem \ref{T51} holds for the model where customers pay at the time of service completion for both infinite-horizon discounted reward and long-run average reward criteria.

\section{Computational results} \label{S6}
In this section, we first provide a policy iteration algorithm that uses our results on the structure of the optimal policy. We then argue that this algorithm with a truncated state space can be used to compute the optimal service rate policy in systems with infinite state space and illustrate this notion with a numerical example. In the interest of space, technical results are only provided for the version of the problem where customers pay at the time of arrival but all results in this section hold for the version where customers pay at time of service completion by modifying the service  cost function as specified in Theorem \ref{T41}. 

Algorithm \ref{Al1} uses the optimization problem in (\ref{NL2}) to compute the long-run average reward optimal policy.
The following result from Guo et al. \cite{survey} provides the convergence properties of Algorithm \ref{Al1}.
\begin{Proposition}
    Let $\{\pi_n\}_{n\geq0}\in\mathbb{F}$ be the sequence of policies obtained by the policy iteration algorithm (Algorithm \ref{Al1}). Then one of the following results holds. Either

    1. the algorithm finitely converges to an optimal long-run average reward policy, or

    2. the sequence $\{g_n\}_{n\geq0}$ converges to the optimal gain $g^*$, and any limit point of $\{\pi_n\}_{n\geq0}$ is long-run average reward optimal.
\end{Proposition}

In Step 8 of Algorithm 6.1, Proposition \ref{PNL2} establishes that we must minimize a convex function. While several convex optimization algorithms guarantee a solution to this problem, the monotonic structure of the optimal policy can be leveraged to further accelerate the search. However, in Step 5, standard computational frameworks cannot directly solve an infinite-dimensional system of equations. To address this tractability issue, we demonstrate that as the state tends to infinity, the optimal service rate converges to the upper bound established in Proposition \ref{PNL2}.

\begin{algorithm}[H]
\caption{Policy Iteration Algorithm}\label{Al1}
\begin{algorithmic}[1]
\STATE Find the function $f^*$ from the projection of the action space $H=\{(\mu(a),f(a)),a\in A\}$.
\STATE Initialize $n\gets1$. Set $\pi_1$ to be a static policy with $\mu(a_i)\gets\mu_-$ for all $i\geq 1$.
\WHILE{True}
    \STATE Compute the stationary distribution $P_{\pi_n}$ under policy $\pi_n$.
    \STATE Solve for vector $u_n$ and constant $g_n$: 
\begin{align*}
   g_n &= r(i,a_i) + \sum_{j\in S} q_{ij}(a_i) u_n(j), && \forall i\in S \\
   0   &= \sum_{i\in S} u_n(i) P_{\pi_n}(i)
\end{align*}
    \STATE Compute $\Delta u_n(i)\gets u_n(i+1)-u_n(i)$

    \STATE Set $i=1$, $A_1=[\mu_-,\mu_+]$.
    \STATE Set 
    \[A_i^*=\mathop{argmax}_{\mu\in A_i} \{-f^*(\mu)-\mu\Delta u_n(i-1)\}.\]
    \STATE Set $A_{i+1}=[\max{\mu\in A^*_i},\mu_+]$.
    \STATE $i\gets i+1$ and return to step 8.
    \STATE Pick $\pi_{n+1}\in$ {\LARGE $\times$}$_{i\geq 1} A_i^*$ (where {\LARGE $\times$} is the Cartesian product), setting $\pi_{n+1}=\pi$ if possible.
    
    \IF {$\pi_n=\pi_{n+1}$}
        \BREAK
    \ELSE
        \STATE $n\gets n+1$
    \ENDIF
\ENDWHILE
\RETURN $(g_n,\ \pi_n)$
\end{algorithmic}
\end{algorithm}

\begin{theorem}\label{T61}
    Let $\pi^*$ be an optimal policy for the long-run average reward criterion and $f^*$ be the lower boundary of the convex hull of the action space. Then, we have
    \[\lim_{i\to\infty}\Delta u(i)=-\frac{h+c\theta}{\theta}\ \mbox{ and }\  \lim_{i\to\infty}\mu(a_i^*)=\sup_\mu\{\partial_{-}f^*(\mu)< \frac{h+c\theta}{\theta}\}.\]
\end{theorem}
\begin{proof}
    Consider the equation in (\ref{EEv3}) for $i>0$. We have
    \begin{equation}\label{OE311S6}
        0=-h-c\theta+\lambda\Delta^2 u(i)-(i-1)\theta\Delta^2 u(i-1)-\theta\Delta u(i)+w^\prime(i).
    \end{equation}
    From Proposition \ref{P35}, we know that $\Delta u(i)$ is bounded below by $-\frac{h+c\theta}{\theta}$. Since $\Delta^2 u(i)\leq 0$ for all $i\in 
    S$, $\Delta u(i)$ is non-increasing in $i$. As a result, we know that $\Delta u(i)$ must converge and $\lim_{i\to\infty}\Delta^2 u(i)=0$. Next, from Lemma \ref{L38C}, we know that $\lim_{i\to\infty }w^\prime(i)=0$. Then, taking the limit of both sides of the equation in  (\ref{OE311S6}), we have
        \[0=-h-c\theta+\lim_{i\to\infty}[\lambda\Delta^2u(i)-(i-1)\theta\Delta^2u(i-1)-\theta\Delta u(i)+w^\prime(i)].\]
Since $\lim_{i\to\infty}\Delta^2 u(i)=0$ and $\lim_{i\to\infty}w^\prime(i)=0$, we know that $(i-1)\theta\Delta^2 u^*(i-1)$ converges as $i$ tends to infinity. Furthermore, the limit should be $0$ because otherwise $\Delta u(i)=\Delta u(0)+\sum_{j=0}^{i-1}\Delta^2 u(j)$ would not converge. As a result, we have
    \[0=-h-c\theta-\theta\lim_{i\to\infty}\Delta u(i)\]
  and hence,
    \[\lim_{i\to\infty}\Delta u(i)=-\frac{h+c\theta}{\theta}.\]
    Since $f^*$ in the convex optimization problem in (\ref{NL2}) is continuous and semi-differentiable, the solution (optimal service rate in state $i\geq 1$) to the optimization problem satisfies
    \[\partial_{-}f^*(\mu)\in \sup\{\partial_{-}f^*(\mu)\leq -\Delta u(i)\}.\]
   By convention, we set $\partial_-f^*(\mu_-)=-\infty$ so that the set $K_i=\{\mu:\partial_{-}f^*(\mu)\in \sup\{\partial_{-}f^*(\mu)\leq -\Delta u(i)\}\}$ is not empty for all $i\geq 1$. Since $\{\Delta u(i)\}_{i\geq0}$ is a non-increasing sequence converging to $-\frac{h+c\theta}{\theta}$, we know that by continuity, the optimal service rate $\{\mu(a_i^*)\}_{i\geq 1}$ is a non-decreasing sequence converging to the infimum of the set 
    \[K_\infty=\{\mu:\sup\{\partial_{-}f^*(\mu)\leq \frac{h+c\theta}{\theta}\}\}.\]
    As a result, we have
   \[ \lim_{i\to\infty}\mu(a_i^*)=\sup_\mu\{\partial_{-}f^*(\mu)< \frac{h+c\theta}{\theta}\}.\]
\end{proof}
For the infinite-horizon discounted criterion, Algorithm \ref{Al1} can be easily modified by changing the system of equations in step 5 into
\[\alpha v_n(i) = r(i,a_i) + \sum_{j\in S} q_{ij}(a_i) v_n(j),\qquad \forall i\in S \\\]
solving for $v_n$, and using the optimization problem in (\ref{NL1}) instead of (\ref{NL2}) in steps 7 and 8. The convergence of the policy iteration algorithm for the infinite-horizon discounted reward criterion follows from Theorem 6.4.6 of Puterman \cite{put}. The following theorem is similar to Theorem \ref{T61}.
\begin{theorem}\label{T62}
    Let $\pi^*$ be an optimal policy under the infinite-horizon discounted reward criterion and $f^*$ be the lower boundary of the convex hull of the action space. Then, we have
    \[\lim_{i\to\infty}\Delta v^*(i)=-\frac{h+c\theta}{\alpha+\theta}\ \mbox{ and }\  \lim_{i\to\infty}\mu(a_i^*)=\sup\{\partial_{-}f^*(\mu)< \frac{h+c\theta}{\alpha+\theta}\}.\]
\end{theorem}
Theorems \ref{T61} and \ref{T62} show that the optimal service rate converges monotonically as the state increases. We then have the following corollary.
\begin{corollary}
    For the infinite-horizon discounted reward criterion, if the lower boundary function $f^*$ is a linear function in $[\mu_-,\mu_+]$ and $\frac{\mathrm{d}f^*}{\mathrm{d}\mu}\geq \frac{h+c\theta}{\alpha+\theta}$, then \[\mu(a_i^*)=\mu_-\quad \forall i\geq 1.\]
    For the long-run average reward criterion,  if the lower boundary function $f^*$ is a linear function in $[\mu_-,\mu_+]$ and $\frac{\mathrm{d}f^*}{\mathrm{d}\mu}\geq \frac{h+c\theta}{\theta}$, then \[\mu(a_i^*)=\mu_-\quad \forall i\geq 1.\]
\end{corollary}

\begin{figure}[H]
\centering
\includegraphics[width=1.05\textwidth,height=0.4\textheight]{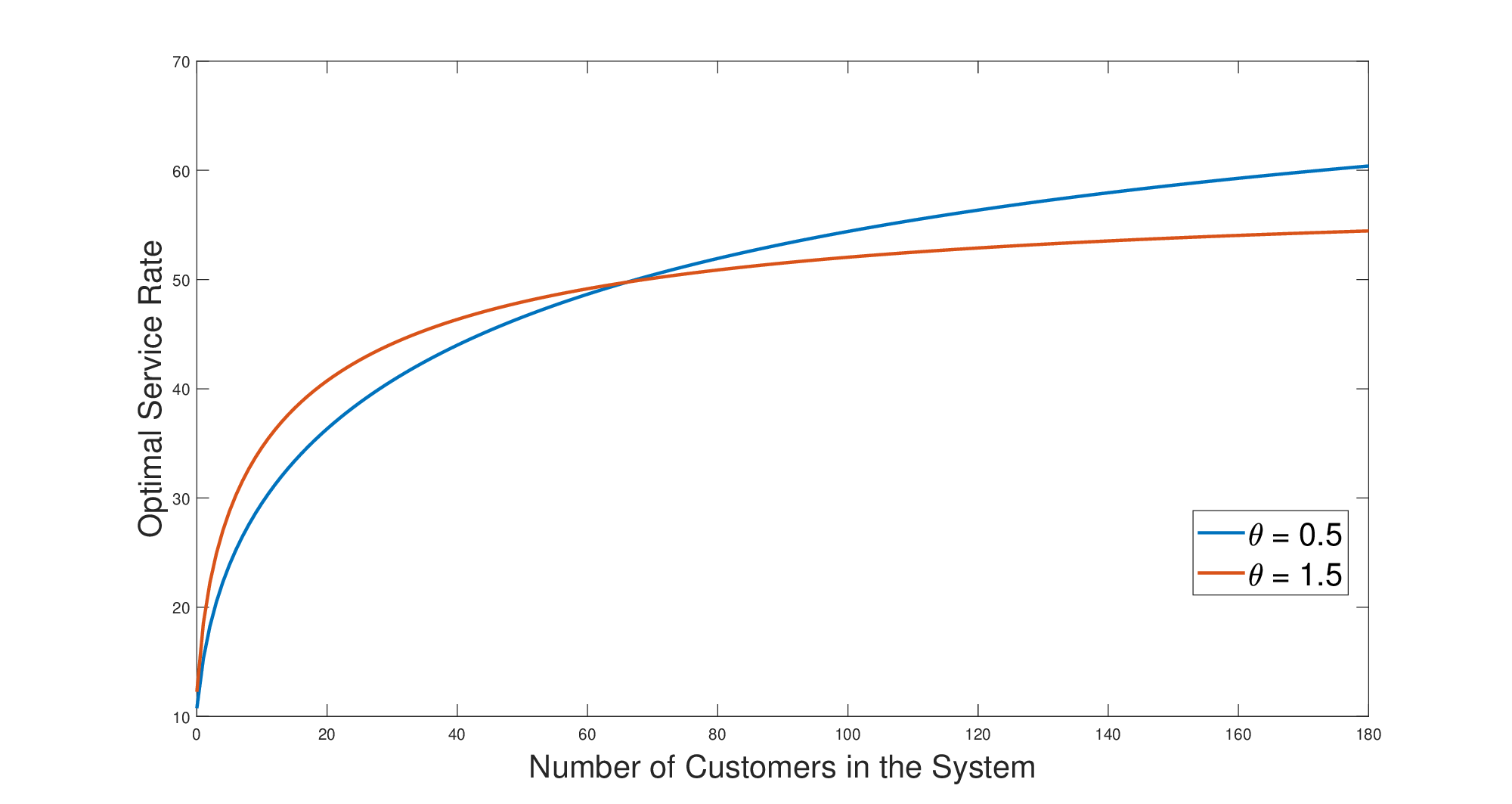}
\caption{Optimal service rate as a function of the number of customers for two values of $\theta$.}
\label{F2}
\end{figure}

\begin{remark}
   If the model does not feature a compact action space, the limits in Theorems \ref{T61} and \ref{T62} may not exist even if the action space is closed. Therefore, compactness is necessary for ensuring the existence of an optimal policy.
\end{remark}
\begin{remark}
    In Theorem \ref{T61}, the term $\frac{h+c\theta}{\theta}$ approaches infinity as $\theta$ tends to $0$. Consequently, under the long-run average reward criterion, the optimal policy for the system without customer abandonments eventually utilizes the maximum possible service rate once the state becomes sufficiently large. In contrast, for systems with abandonments, the optimal service rate converges to a finite value as the state tends to infinity. Furthermore, this asymptotic limit is monotonically non-increasing in the abandonment rate $\theta$, implying that the optimal service rate in any given state is not necessarily non-decreasing in $\theta$; for sufficiently large states, it may instead decrease as the abandonment rate increases. For example, consider a system with $\lambda=10,\ h=2,\ c=8,\ \mathcal{A}=[1,100],\ \mu(a)=a$, and $f(a)=0.08a^2$. Figure \ref{F2} depicts the optimal service rate as a function of the number of customers for this system for two values of the abandonment rate: $\theta=0.5$ and $\theta=1.5$. The figure illustrates that a higher abandonment rate does not necessarily result in a higher optimal service rate.
\end{remark}

\color{black}

\begin{figure}[H]
\centering
\includegraphics[width=1.05\textwidth,height=0.4\textheight]{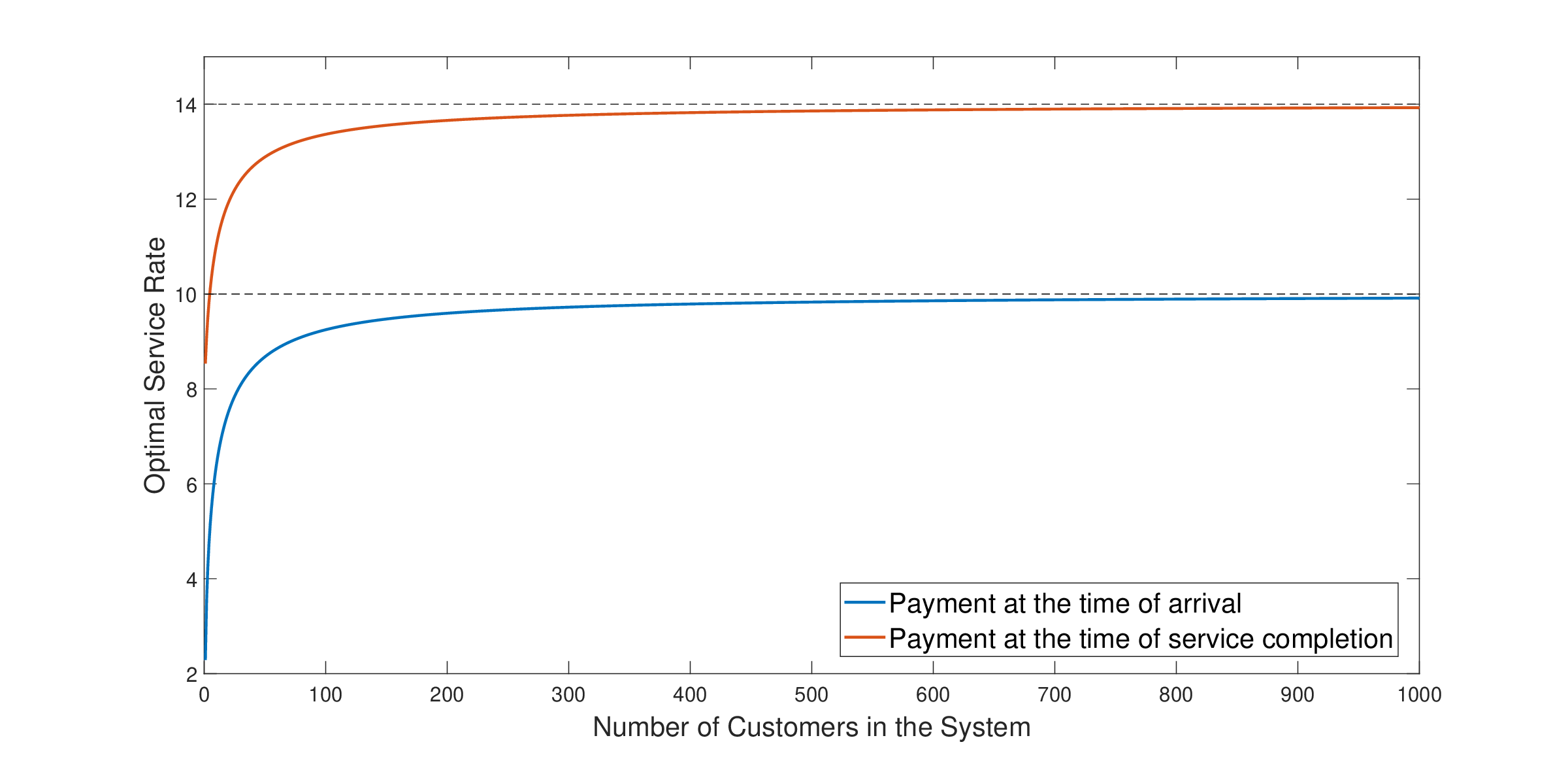}
\caption{Optimal service rate in each state for both versions of the problem.}
\label{F1}
\end{figure}

We conclude this section with a numerical example for which we consider both versions of the problem (payment at the time of arrival and service completion) and compute the  long-run average reward optimal policy. Note that when implementing Algorithm \ref{Al1} for a system with an infinite state space, one needs to truncate the state space. While we do not provide a formal truncation bound at this stage and leave this important question as a topic for future research, we suggest choosing a sufficiently large truncation level. If the policy obtained from the policy iteration algorithm is close to a monotone policy—meaning the service rate converges to the results in Theorem \ref{T62}—then we can conclude that the chosen truncation level is appropriate. Consider a system with $\lambda=0.5$, $r=2$, $h=1$, $c=3$, $\theta=0.5$, $\mathcal{A}=[0.5,30]$, $\mu(a)=a$ and $f(a)=0.25a^2$. We use Algorithm \ref{Al1} with a truncated state space $\{0,\dots,1000\}$ to compute the long-run average reward optimal policy. Figure \ref{F1} depicts the optimal service rate as a function of the state for both versions.

As we have shown, Figure \ref{F1} illustrates that optimal service rate is monotone non-decreasing in the state and converges to the limit specified in Theorem \ref{T61}. Furthermore, the optimal service rate in the version with payment after service completion is higher than the one in the version with payment at the time of arrival (see Corollary \ref{C45}).

\section{Conclusion} \label{S7}
This paper studies the service rate control problem of a Markovian queueing system with abandonments. Apart from compactness and continuity, no additional assumptions are imposed on the action space.
Under the assumptions that there is a holding cost and a cost associated with abandonments and a reward obtained either at the time of arrivals or service completions, we compute the optimal service rate policy under both infinite-horizon discounted and long-run average reward optimality criteria. When the buffer size is infinite, we formulate the queueing system as a Continuous Time Markov Decision Process (CTMDP) (because the problem is not uniformizable). From the optimality equations of CTMDP, using limit theorems, we directly show the monotonicity of service rate with respect to the number of customers in the system. Furthermore, we prove that the optimal service rates lie on the lower boundary of the convex hull of the action space.

We also consider systems with finite buffer and show the interesting result that optimal service rate is not necessarily monotone in the state. Specifically, we show that the optimal service rate first increases and then decreases with respect to the number of customers in the system. Furthermore, we prove that the optimal service rate converges as the state goes to infinity. This property allows us to truncate the state space at suitable levels in order to numerically compute the optimal service rate in infinite buffer models.

As a potential avenue for future work, one could consider the joint control of admission and service rates in systems with customer abandonments. While we anticipate that our solution techniques are applicable to this extended environment, the structural properties of the optimal policy remain to be fully characterized. A double-sets -of-thresholds structure was established by Wu and Ayhan \cite{Adm} for the pure admission control problem; however, whether the co-existence of service rate control preserves this threshold structure or yields a more complex policy is a non-trivial open problem.

\section{Acknowledgement}
The work of the second author was supported by the National Science Foundation Grant CMMI-2127778.

\section*{Appendix}
\noindent The next lemma provides another equivalency between the two versions of the problem under the long-run average reward optimality.

\begin{Alemma}
    For the long-run average reward criterion, the version with customers paying reward $r$ at the time of service completion is equivalent to the version with customers paying reward $r$ at the time of arrival and abandonment cost $c^\prime=c+r$.
\end{Alemma}
\begin{proof}
    From Theorem \ref{A4} (i) we know that the under any policy $\pi\in\mathbb{F}$ the continuous time Markov chain has a unique invariant probability measure $P_\pi=(P_0,P_1,\dots)$. Then, the long-run average reward $\mathcal{G}(i,\pi)$ can be written as
    \[\mathcal{G}=\sum_{i=0}^\infty r^\prime(i,a_i)P_i.\]
    Note that under policy $\pi$, the two models (corresponding to paying at the time of arrivals and at the time of service completion) have the same transition rate matrix and, hence, the same invariant probability measure $P_\pi$. Therefore, in the model with customers paying at the time of arrival with abandonment cost $c^\prime=c+r$, we have 
    \begin{eqnarray}
        g&=&\sum_{i=0}^\infty r(i,a_i)P_i\nonumber\\
        &=&P_0\lambda r+\sum_{i=1}^\infty P_i(\lambda r-hi-c^\prime(i-1)\theta-f(a_i))\nonumber\\
        &=&P_0\lambda r+\sum_{i=1}^\infty P_i(\lambda r-hi-(c+r)(i-1)\theta-f(a_i))\nonumber\\
        &=&\sum_{i=0}^\infty P_i\lambda r-\sum_{i=1}^\infty P_i(i-1)\theta r-\sum_{i=1}^\infty P_i(hi+c(i-1)\theta+f(a_i)).\nonumber
    \end{eqnarray}
    Note that $P_0\lambda=P_1\mu(a_1)$ and for $i\geq1$, $P_{i}\lambda=P_{i+1}(\mu(a_{i+1})+(i-1)\theta)$. Therefore, 
    \[\sum_{i=0}^\infty P_i\lambda r=P_1\mu(a_1)+\sum_{i=1}^\infty P_i(\mu(a_{i})+(i-1)\theta).\]
    As a result,
    \begin{eqnarray}
    g&=&\sum_{i=1}^\infty P_i\mu(a_i) r-\sum_{i=1}^\infty P_i(hi+c(i-1)\theta+f(a_i))\nonumber\\
    &=&\sum_{i=1}^\infty P_i(\mu(a_i) r-hi-c(i-1)\theta-f(a_i))\nonumber\\
    &=&\sum_{i=1}^\infty P_ir^\prime(i,a_i)=\mathcal{G}\nonumber
     \end{eqnarray}
     which completes the proof.
\end{proof}

\begin{thebibliography}{99}
\bibitem{john} Adusumilli, K. M. and Hasenbein, J. J. (2013) ``Dynamic Admission and Service Rate Control of a Queue,'' {\em Queueing Systems}, 66, 131--154.
\bibitem{mor} Armony, M., Chan, C. W., and Zhu, B. (2017) ``Critical Care Capacity Management: Understanding the Role of a Step Down Unit,'' {\em Production and Operations Management}, 27, 859--883.
\bibitem{baris} Ata, B. and Shneorson, S. (2006) ``Dynamic Control of an M/M/1 Service System with Adjustable Arrival and Service Rates,'' {\em Management Science}, 52, 1778--1791.
\bibitem{Atar} Atar R., Giat C. and Shimkin N. (2010). "The $c\mu/\theta$ rule for many-server queues with abandonment," {\em Operations Research}, 58, 1427-1439.
\bibitem{pamela} Badian-Pessot, P., Lewis, M. E., and Down, D. G. (2021) ``Optimal Control Policies for an M/M/1 Queue with a Removable Server and Dynamic Service Rates,'' {\em Probability in the Engineering and  Informational Sciences} 35, 189--209.
\bibitem{batt} Batt, R. J. and Terwiesch, C. (2015) ``Waiting Patiently: An Empirical Study of Queue Abandonment in an Emergency Department,'' {\em Management Science}, 61, 39--59.

\bibitem{sunjay} Bhulai, S., Blok, H., and Spieksma, F. M. (2022) ``$K$ competing queues with customer abandonment: optimality of a generalised $c \mu$- rule by the Smoothed Rate Truncation method'', {\em Annals of Operations Research}, 317, 387--416.

\bibitem{bora} \c{C}ekyay, B. (2024) ``Discounted Cost Exponential Semi-Markov Decision Processes with Unbounded Transition Rates: a Service Rate Control Problem with Impatient Customers,'' {\em Probability in the Engineering and Informational Sciences}, 38(4), 668--694. 
\bibitem{chen} Chen, G.,  Liu, Z., and Xia, L. (2023) ``Event-based Optimization of Service Rate Control in Retrial Queues,'' {\em Journal of the Operational Research Society}, 74, 979--991.
\bibitem{crab1} Crabill, T. B., (1972) ``Optimal Control of a Service Facility with
Variable Exponential Service Times and Constant Arrival Rate,'' {\em Management Science}, 18, 560--566.
\bibitem{crab2} Crabill, T. B. (1974) ``Optimal Control of a Maintenance System with Variable Service Rates,'' {\em Operations Research}, 22, 736--745.

\bibitem{Down} Down, D. G., Koole, G. and Lewis, M. E. (2011) ``Dynamic control of a single-server system with abandonments," {\em Queueing Systems}, 67, 63–90.

\bibitem{george} George, J. M. and  Harrison, M. J. (2001) ``Dynamic Control of a Queue with Adjustable Service Rate," {\em Operations Research}, 49, 720--731.
\bibitem{gong} Gong, J. and Li, M. (2014) ``Queueing Time Decision Model with the Consideration on Call Center Customer Abandonment Behavior,'' {\em Journal of Networks}, 9, 2441--2447.
\bibitem{survey} Guo, X., Hernández-Lerma, O., and Prieto-Rumeau, T. (2006). ``A survey of recent results on continuous-time Markov decision processes," {\em TOP 14, 177–261}.

\bibitem{Kocaga} Ko\c{c}a\u{g}a, Y. L. (2017). ``An approximating diffusion control problem for dynamic admission and service rate control in a $G/M/N+G$ queue," {\em Operations Research Letters}, 45(6), 538-542.

\bibitem{ger} Koole, G. and Pot, A. (2011) ``A Note on Profit Maximization and Monotonicity for Inbound Call Centers,'' {\em Operations Research}, 59, 1304--1308.
\bibitem{jo} Jo, K. Y. and Stidham, S. (1983) ``Optimal Service-rate Control of M/G/1 Queueing Systems Using Phase Methods,'' {\em Advances in Applied Probability}, 15, 616--637.
\bibitem{ravi} Kumar, R.,  Lewis, M. E., and   Topalo\u{g}lu, H. (2013) ``Dynamic Service Rate Control for a Single-server Queue with Markov-modulated Arrivals,'' {\em Naval Research Logistics}, 60, 661--677.
\bibitem{kulkarni} Nelson, L. and Kulkarni, V. G. (2014) ``Optimal Arrival Rate and Service Rate Control of Multi-server Queues,'' {\em Queueing Systems}, 76, 37--50.
\bibitem{put} Puterman, M.~L. {\em Markov Decision Processes}\/. John Wiley \& Sons, New York, NY, 1994.
\bibitem{convex} Rockafellar, R. T., \textit{Convex analysis}, Princeton University Press, 1970.
\bibitem{rudin}  Rudin, W., \textit{Principles of Mathematical Analysis}, 3rd ed., McGraw-Hill, 1976.
\bibitem{sarhangian} Sarhangian, V., Abouee-Mehrizi, H., Baron, O. and Berman, O. (2018) ``Threshold-Based Allocation Policies for Inventory Management of Red Blood Cells,'' {\em Manufacturing \& Service Operations Management}, 20, 347--362.
\bibitem{dick} Serfozo, R. F. (1981) ``Optimal Control of Random Walks, Birth and Death Processes, and Queues,'' {\em Advances in Applied Probability}, 13, 61--83.
\bibitem{deniz} \c{S}im\c{s}ek, D., Wu, C., Bassamboo, A., and Perry, O. (2025) ``A Unified Fluid Model for Large Service Systems with Patience- or Delay-Dependent Service Times,'' {\em Available at SSRN}.
\bibitem{stidham} Stidham, S. and   Weber, R. R. (1989) ``Monotonic and Insensitive Optimal Policies for Control of Queues with Undiscounted Costs,'' {\em Operations Research}, 37, 611--625.
\bibitem{weber} Weber, R. R. and   Stidham, S. (1987) ``Optimal Control of Service
Rates in Networks of Queues,'' {\em Advances in Applied Probability}, 19, 202--218.
\bibitem{wu} Wu, H., He, Q. and Erenay, F. S. (2025) ``Double‐sided Queues and Their Applications to Vaccine Inventory Management,'' {\em Naval Research Logistics}, 72, 292--316.

\bibitem{Adm} Wu R. and  Ayhan, H. (2025) ``Optimal admission control in queues with multiple customer classes and abandonments,'' Queueing Systems, 109(1): 6.

\bibitem{gordon} Zheng, Y.,  Julaiti, J., and Pang, G. (2024) ``Adaptive Service Rate Control of an M/M/1 Queue with Server Breakdowns,'' {\em Queueing Systems}, 106, 159--191.




\end{thebibliography}
\end{document}